\documentclass[11pt]{article}
\usepackage{cite}
\usepackage{mathrsfs}
\usepackage{amsfonts}
\usepackage{amsmath}
\usepackage{amsfonts,amssymb,color}
\usepackage{dsfont}
\usepackage{curves}
\usepackage{mathrsfs}
\usepackage{pifont}
\usepackage{amssymb}
\usepackage{latexsym,amsmath,amssymb,amsfonts,epsfig,graphicx,cite,psfrag}
\usepackage{eepic,color,colordvi,amscd}

\newtheorem{theorem}{Theorem}[section]

\newtheorem{lemma}{Lemma}[section]

\newtheorem{claim}{Claim}[section]
\newtheorem{conjecture}{Conjecture}[section]

\newcommand{\qed}{\hfill\rule{0.5em}{0.809em}}

\def\emptyset{\mbox{{\rm \O}}}

\renewcommand{\baselinestretch}{1.1}
\def\qed{\hfill \rule{4pt}{7pt}}

\def\pf{\noindent {\it Proof. }}

\begin{document}

\title{The chromatic number of heptagraphs\thanks{Partially supported by NSFC projects 11931006 and 12101117, and NSFJS No. BK20200344}}
 \author{Di Wu$^{1,}$\footnote{Email: 1975335772@qq.com},  \;\; Baogang  Xu$^{1,}$\footnote{Email: baogxu@njnu.edu.cn OR baogxu@hotmail.com.}\;\; and \;\; Yian Xu$^{2,}$\footnote{Email: yian$\_$xu@seu.edu.cn. }\\\\
\small $^1$Institute of Mathematics, School of Mathematical Sciences\\
\small Nanjing Normal University, 1 Wenyuan Road,  Nanjing, 210023,  China\\
\small $^2$School of Mathematics, Southeast University, 2 SEU Road, Nanjing, 211189, China}
\date{}

\maketitle


\begin{abstract}
A hole is an induced cycle of length at least 4. A graph is called a pentagraph if it has no cycles of length 3 or 4 and has no holes of odd length at least 7, and is called a heptagraph if it has no cycles of length less than 7 and has no holes of odd length at least 9. Let $\l\ge 2$ be an integer. The current authors proved that a graph is 4-colorable if it has no cycles of length less than $2\l+1$ and has no holes of odd length at least $2\l+3$. Confirming a conjecture of Plummer and Zha, Chudnovsky and Seymour proved that every pentagraph is 3-colorable. Following their idea, we show that every heptagraph is 3-colorable.
\begin{flushleft}
\end{flushleft}

\end{abstract}

\newpage

\section{Introduction}

Let $G$ be a graph, and let $u$ and $v$ be two vertices of $G$. We simply write $u\sim v$ if $uv\in E(G)$, and write $u\not\sim v$ if $uv\not\in E(G)$. We use $d_G(u)$ (or simply $d(u)$) to denote the degree of $u$ in $G$, and let $\delta(G)=\min\{d(u):u\in V(G)\}$. Let $S$ be a subset of $V(G)$. We use $G[S]$ to denote the subgraph of $G$ induced by $S$.  For two graphs $G$ and $H$, we say that $G$ induces $H$ if $H$ is an induced subgraph of $G$.

Let $S$ and $T$ be two subsets of $V(G)$, and let $x$ and $y$ be two vertices of $G$.  We use $N_S(x)$ to denote the neighbors of $x$ in $S$, and define $N_S(T)=\cup_{x\in T} N_S(x)$ (if $S=V(G)$ then we omit the subindex and simply write them as $N(x)$ or $N(T)$). An $xy$-path is a path between $x$ to $y$, and an $(S, T)$-path is a path $P$ with  $|S\cap P|=|T\cap P|=1$.   A cycle on $k$ vertices is denoted by $C_k$. Let $P$ be a path, we use $\l(P)$ and $P^*$ to denote the length and the set of internal vertices of $P$, respectively. If $u, v\in V(P)$, then  $P[u, v]$ denotes the segment of $P$ between $x$ and $y$.

Let $k$ be a positive integer.  A {\em hole}  is an induced cycle of length at least 4, a hole of length $k$ is called a $k$-hole, and a $k$-hole is said to be an {\em odd} (resp. {\em even})  hole if $k$ is odd (resp. even).

A $k$-{\em coloring} of $G$ is a mapping $c: V(G)\mapsto \{1, 2, \ldots, k\}$ such that $c(u)\neq c(v)$ whenever $u\sim v$ in $G$.  The {\em chromatic number} $\chi(G)$ of $G$ is the minimum integer $k$ such that $G$ admits a $k$-coloring.

Let $\l\ge 2$ be an integer. Let ${\cal G}_l$ denote the family of graphs which have no cycles of length less than $2\l+1$ and have no odd holes of length at least $2\l+3$. The graphs in ${\cal G}_2$ are called {\em pentagraphs}, and the graphs in ${\cal G}_3$ are called {\em heptagraphs}.

A 3-connected graph is said to be {\em internally 4-connected} if every cutset of size 3 is the neighbor set of a vertex of degree 3. Robertson conjectured (see \cite{NPRZ2011}) that the Petersen graph is the only non-bipartite pentagraph which is 3-connected and internally 4-connected. In 2014, Plummer and Zha \cite{MPXZ} presented some counterexamples to Robertson's conjecture, and posed the following new conjecture.

\begin{conjecture}\label{conj-P-Z}{\em (\cite{MPXZ})}
Every pentagraph is $3$-colorable.
\end{conjecture}

Xu, Yu and Zha \cite{XYZ2017} proved that every pentagraph is 4-colorable. Very recently,  Chudnovsky and Seymour \cite{MCPS2022} presented a structural characterization for pentagraphs and confirmed Conjecture~\ref{conj-P-Z}.

A $P_3$-{\em cutset} of $G$ is an induced path $P$ on three vertices such that $V(P)$ is a cutset. A {\em parity star}-{\em cutset} is a cutset $X\subseteq V(G)$ such that $X$ has a vertex, say $x$, which is adjacent to every other vertex in $X$, and $G-X$ has a component, say  $A$, such that every two vertices in $X\setminus\{x\}$ are joint by an induced even path with interior in $V(A)$.

\renewcommand{\baselinestretch}{1}
\begin{theorem}\label{theo-1-1}{\em (\cite{MCPS2022})}
Let $G$ be a pentagraph which is not the Petersen graph. If $\delta(G)\ge 3$, then $G$ is either bipartite, or admits a $P_3$-cutset or a parity star-cutset.
\end{theorem}\renewcommand{\baselinestretch}{1.2}

Below is a lemma contained in the proof of  \cite[Theorem 1.1]{MCPS2022}.

\renewcommand{\baselinestretch}{1}
\begin{lemma}\label{lem-critical} {\em (\cite{MCPS2022})}
Let $G$ be a pentagraph that is not the Petersen graph. If $\chi(G)=4$ and every proper induced subgraph of $G$ is $3$-colorable, then $G$ has neither $P_3$-cutsets nor parity-star cutsets.
\end{lemma}\renewcommand{\baselinestretch}{1.2}

As a direct consequence of Theorem~\ref{theo-1-1} and Lemma~\ref{lem-critical}, one can easily verify that every pentagraph is 3-colorable. Thus, Conjecture~\ref{conj-P-Z} is true.

By generalizing the result of \cite{XYZ2017},  the current authors \cite{WXX2022} proved that $\chi(G)\le 4$ for each graph $G\in \cup_{\l\ge 2}{\cal G}_l$, and conjectured that $\chi(G)\le 3$ for such graphs.

\begin{theorem}\label{theo-1-2-0}{\em (\cite{WXX2022})}
All graphs in $\cup_{\l\ge 2}{\cal G}_l$ are  $4$-colorable.
\end{theorem}

In this paper, we prove the conjecture for heptagraphs.

\begin{theorem}\label{theo-1-2}
Every heptagraph is $3$-colorable.
\end{theorem}

We follow the idea of Chudnovsky and Seymour and prove  a structural theorem for heptagraphs.

\renewcommand{\baselinestretch}{1}
\begin{theorem}\label{theo-1-3}
Let $G$ be a heptagraph. If $\delta(G)\ge 3$, then $G$ is bipartite, or
admits a $P_3$-cutset or a parity star-cutset.
\end{theorem}\renewcommand{\baselinestretch}{1.2}

A  conclusion similar to Lemma~\ref{lem-critical} also holds on graphs in $\cup_{\l\ge 2}{\cal G}_l$. Since its proof is almost the same as that of Lemma~\ref{lem-critical}, we leave the proof to readers.

\renewcommand{\baselinestretch}{1}
\begin{lemma}\label{lem-critical-H}
Let $G$ be a graph in $\cup_{\l\ge 2}{\cal G}_l$. If $\chi(G)=4$ and every proper induced subgraph of $G$ is $3$-colorable, then $G$ has neither $P_3$-cutsets nor parity-star cutsets.
\end{lemma}\renewcommand{\baselinestretch}{1.2}

\noindent{\bf Assuming Theorem~\ref{theo-1-3}, we can prove Theorem~\ref{theo-1-2}}: Suppose to its contrary, let $G$ be a heptagraph with $\chi(G)=4$ such that all proper induced subgraphs of $G$ are 3-colorable. It is certain that $G$ is not bipartite, and $\delta(G)\ge 3$. By Theorem~\ref{theo-1-3}, we have that $G$ must have a $P_3$-cutset or a parity-star cutset, which leads to a contradiction to Lemma~\ref{lem-critical-H}. Therefore, Theorem~\ref{theo-1-2} holds. \qed

In Sections 2 and 3, we discuss the structure of heptagraphs and prove some lemmas. Theorem~\ref{theo-1-3} is proved in Section 4.

\section{Subgraphs ${\cal P}$ and ${\cal P}'$}

If a cutset is a clique, then we call it a {\em clique cutset}. It is certain that in any triangle free graph, every clique cutset is a parity star-cutset.

Let $H$ be a proper induced subgraph of $G$ and $s, t\in V(H)$ with $s\not\sim t$, let $P$ be an induced $st$-path such that $\l(P)\ge 3$ and $P^*\subseteq V(G)\setminus V(H)$. If every vertex of $H-\{s, t\}$ that has a neighbor in $P^*$ is adjacent to both $s$ and $t$, then we call $P$ an $st$-{\em ear} of $H$.
\renewcommand{\baselinestretch}{1}
\begin{lemma}\label{MCPS-lem-2-1}{\em (\cite[2.1]{MCPS2022})}
Let $G$ be a pentagraph without clique cutsets, and let $H$ be a proper induced subgraph of $G$ with $|V(H)|\ge3$. If each vertex of $V(G)\setminus V(H)$ has at most one neighbor in $V(H)$, then there exist nonadjacent vertices $s,t\in V(H)$ such that $H$ has an $st$-ear.
\end{lemma}\renewcommand{\baselinestretch}{1.2}

A similar conclusion holds to heptagraphs. We omit its proof as it is the same as that proof of Lemma~\ref{MCPS-lem-2-1} of \cite{MCPS2022}.

\renewcommand{\baselinestretch}{1}
\begin{lemma}\label{lem-2-1}
Let $G$ be a heptagraph without clique cutsets, and let $H$ be a proper induced subgraph of $G$ with $|V(H)|\ge3$. If each vertex of $V(G)\setminus V(H)$ has at most one neighbor in $V(H)$, then there exist nonadjacent vertices $s,t\in V(H)$ such that $H$ has an $st$-ear.
\end{lemma}\renewcommand{\baselinestretch}{1.2}

A {\em big odd hole} is an odd hole of length at least 9, and a {\em short cycle} is a cycle of length at most 6.

\begin{figure}[htbp]\label{fig-4}
	\begin{center}
		\includegraphics[width=8cm]{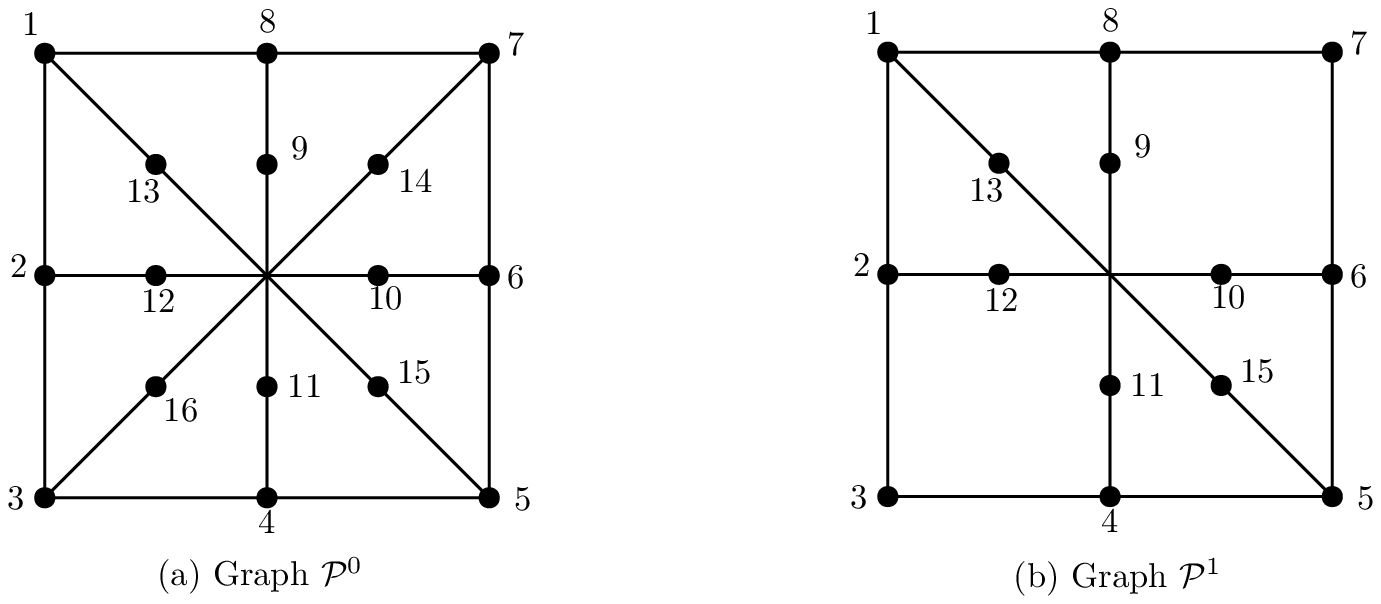}
	\end{center}
	\vskip -15pt
	\caption{Graphs ${\cal P}^0$ and ${\cal P}^1$}
\end{figure}

\begin{figure}[htbp]\label{fig-3}
	\begin{center}
		\includegraphics[width=8cm]{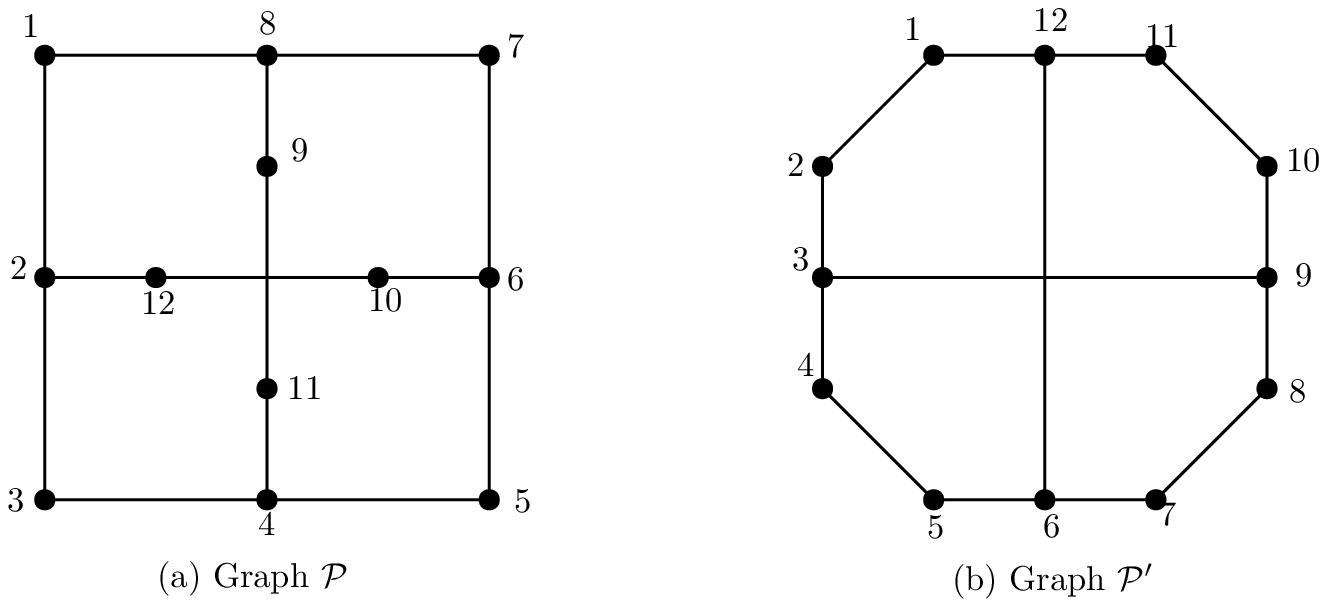}
	\end{center}
	\vskip -15pt
	\caption{Graphs ${\cal P}$ and ${\cal P}'$}
\end{figure}

Let $C_8=1-2-3-4-5-6-7-8-1$, and $C_{12}=1-2-3-4-5-6-7-8-9-10-11-12$. We use ${\cal P}^0$ to denote the graph obtained from $C_8$ by adding four induced paths 1-13-15-5, 2-12-10-6, 3-16-14-7, and 4-11-9-8 (see Figure~\ref{fig-4}$(a)$), use ${\cal P}^1$ to denote the graph obtained from ${\cal P}^0$ by deleting vertices 14 and 16 (see Figure~\ref{fig-4}$(b)$), use ${\cal P}$ to denote the graph obtained from ${\cal P}^1$ by deleting vertices 13 and 15 (see Figure~\ref{fig-3}$(a)$), and use ${\cal P}'$ to denote the graph obtained from $C_{12}$ by adding edges 3-9 and 6-12 (see Figure~\ref{fig-3}$(b)$).

Let $G$ be a non-bipartite heptagraph with $\delta(G)\ge 3$ and without $P_3$-cutsets or parity star-cutsets. We will show that $G$ does not induce ${\cal P}$ or ${\cal P}'$. Firstly, we prove that $G$ does not induce ${\cal P}^0$ or ${\cal P}^1$.

\begin{lemma}\label{theo-2-1}
Let $G$ be a heptagraph that induces a ${\cal P}^0$. If $G$ has no clique cutsets, then $G={\cal P}^0$.
\end{lemma}
\pf Let $H$ be an induced subgraph of $G$ isomorphic to ${\cal P}^0$. Since the distance of any two vertices of $H$ is at most four, no vertex in $V(G)\setminus V(H)$ has more than one neighbor in $V(H)$. We may suppose that $G$ has no clique cutsets and $G\ne H$. By Lemma~\ref{lem-2-1}, there are nonadjacent vertices $s, t\in V(H)$ and an $st$-ear $P$ in $H$. We choose $s$ and $t$ that minimize $\l(P)$. By symmetry, we only need to verify that $P$ is an $st$-ear for $s\in\{1, 8, 9, 13\}$.

If $(s, t)=(1,7)$, let $C=1-P-7-6-10-12-2-1$ and $C'=1-P-7-6-5-4-3-2-1$.  If $(s, t)=(1, 6)$, let $C=1-P-6-10-12-2-1$ and $C'=1-P-6-5-4-3-2-1$. If $(s, t)=(1,5)$, let $C=1-P-5-6-10-12-2-1$ and $C'=1-P-5-4-3-2-1$. If $(s, t)=(1, 9)$, let $C=1-P-9-11-4-3-2-1$ and $C'=1-P-9-11-4-5-6-10-12-2-1$. If $(s, t)=(1, 14)$, let $C=1-P-14-16-3-2-1$ and $C'=1-P-14-16-3-4-11-9-8-1$. If $(s, t)=(1, 10)$, let $C=1-P-10-6-7-8-1$ and $C'=1-P-10-6-5-4-11-9-8-1$. If $(s, t)=(1, 15)$, let $C=1-P-15-5-4-3-2-1$ and $C'=1-P-15-5-4-11-9-8-1$. One can check that either $C$ or $C'$ is a big odd hole in all cases.

Thus we suppose by symmetry that $\{s, t\}\cap \{1, 3, 5 ,7\}=\emptyset$ .

If $(s, t)=(8, 6)$, let $C=8-P-6-10-12-2-1-8$ and $C'=8-P-6-10-12-2-3-4-11-9-8$. If $(s, t)=(8, 4)$, let $C=8-P-4-5-6-7-8$ and $C'=8-P-4-3-2-12-10-6-7-8$. If $(s, t)=(8, 13)$, let $C=8-P-13-15-5-6-7-8$ and $C'=8-P-13-15-5-4-11-9-8$. If $(s, t)=(8, 12)$, let $C=8-P-12-10-6-7-8$ and $C'=8-P-12-10-6-5-4-11-9-8$. If $(s, t)=(8, 16)$, let $C=8-P-16-3-2-1-8$ and $C'=8-P-16-3-4-11-9-8$. If $(s, t)=(8, 11)$, let $C=8-P-11-4-3-2-1-8$ and $C'=8-P-11-4-5-6-10-12-2-1-8$. In all cases, one of $C$ and $C'$ is a big odd hole.

We may further suppose, by symmetry, that $\{s, t\}\cap \{1, 2, 3, 4, 5, 6, 7, 8\}=\emptyset$.

If $(s, t)=(13, 9)$, let $C=13-P-9-11-4-5-15-13$ and $C'=13-P-9-8-7-6-5-15-13$. If $(s, t)=(13, 14)$, let $C=13-P-14-16-3-2-1-13$ and $C'=13-P-14-12-3-4-5-15-13$. If $(s, t)=(13, 10)$, let $C=13-P-10-6-7-8-1-13$ and $C'=13-P-10-6-5-4-11-9-8-1-13$. In each case, one of $C$ and $C'$ must be a big odd hole.

By symmetry, it remains to consider that $(s, t)=(9, 12)$. Now, $l(P)\ge3$,  and either $9-P-12-10-6-7-8-9$ or $9-P-12-10-6-5-4-11-9$ is a big odd hole. This proves Lemma~\ref{theo-2-1}.  \qed

Notice that ${\cal P}^1$ is an induced subgraph of ${\cal P}^0$.
With almost the same arguments, one can check that the following lemma holds.

\renewcommand{\baselinestretch}{1}
\begin{lemma}\label{theo-2-2}
Let $G$ be a heptagraph with no clique cutsets. If $G$ induces a ${\cal P}^1$ then $G\in \{{\cal P}^0, {\cal P}^1\}$.
\end{lemma}	\renewcommand{\baselinestretch}{1.2}

\renewcommand{\baselinestretch}{1}
\begin{lemma}\label{theo-1-4}
Let $G$ be a heptagraph with no clique cutsets. If $G$ induces a ${\cal P}$ then $\delta(G)\le 2$.
\end{lemma}\renewcommand{\baselinestretch}{1.2}
\pf Let $H$ be an induced subgraph of $G$ isomorphic to ${\cal P}$. Suppose that $\delta(G)\ge 3$ and $G$ has no clique cutsets. Since the distance of any two vertices of $H$ is at most four, no vertex in $V(G)\setminus V(H)$ has more than one neighbor in $V(H)$.  By Lemma~\ref{lem-2-1}, there exist nonadjacent vertices $s,t\in V(H)$ such that $H$ has an $st$-ear.

If $(s, t)\in\{(1, 5), (3, 7)\}$ and $\l(P)=3$, then $G[V(H)\cup V(P)]\in \{{\cal P}^0,  {\cal P}^1\}$. By Lemmas~\ref{theo-2-1} and \ref{theo-2-2}, we have that $G\in \{{\cal P}^0,  {\cal P}^1\}$, a contradiction. Suppose that either $(s, t)\not\in \{(1, 5), (3, 7)\}$ or $\l(P)>3$. Notice that ${\cal P}$ is an induced subgraph of ${\cal P}^0$. With totally the same arguments as that used in the proof of Lemma~\ref{theo-2-1}, we can always find a big odd hole in $G$. Therefore, Lemma~\ref{theo-1-4} holds. \qed

\renewcommand{\baselinestretch}{1}
\begin{lemma}\label{theo-1-5}
Let $G$ be a heptagraph that induces a ${\cal P}'$. If $\delta(G)\ge 3$ then $G$ admits a clique cutset or a $P_3$-cutset.
\end{lemma}	\renewcommand{\baselinestretch}{1.2}
\pf Let $H$ be an induced subgraph of $G$ isomorphic to ${\cal P}'$. Since the distance of any two vertices of $H$ is at most four, no vertex in $V(G)\setminus V(H)$ has more than one neighbor in $V(H)$. We assume $\delta(G)\ge 3$  and $G$ does not admit a clique cutset or a $P_3$-cutset. It is certain that $G\ne H$. By Lemma~\ref{theo-1-4}, we may assume that $G$ induces no ${\cal P}$.

Let us call the four sets $\{2,3,4\}, \{5,6,7\}, \{8,9,10\}$, and $\{1,11,12\}$ the {\em sides} of $H$. Since $G$ has no $P_3$-cutsets, there is a connected subgraph $F$ of $G - V(H)$ such that $N_H(F)$ is not a subset of any side of $H$. Choose an $F$ with $|V(F)|$ minimal. Since $N_H(F)$ cannot be a clique, there exist nonadjacent vertices $s, t\in N_H(F)$ and an  induced $st$-path $P$ with $P^*\subseteq V(F)$. We choose $P$ that minimizes $\l(P)$, then every vertex of $H-\{s, t\}$ with a neighbor in $P^*$ is adjacent to both $s$ and $t$.

By symmetry, we may assume that $s\in\{1, 12\}$.

Firstly, suppose that $s=1$ and $t\ne 11$.
If $(s, t)=(1, 8)$ and $l(P)=3$, then $V(P)\cup\{2,3,4,5,6,7,9,12\}$ induces a ${\cal P}$ in $G$, contradicting Lemma~\ref{theo-1-4}. If $(s, t)=(1, 8)$ and $l(P)>3$ then  either $1-P-8-9-10-11-12-1$ or $1-P-8-9-3-2-1$ is a big odd hole, contradicting the choice of $G$. Thus, we assume that $t\ne 8$. If $(s, t)=(1, 10)$, let $C=1-P-10-9-3-2-1$ and $C'=1-P-10-9-3-4-5-6-12-1$. If $(s, t)=(1, 9)$, let $C=1-P-9-10-11-12-1$ and $C'=1-P-9-8-7-6-12-1$.  If $(s, t)=(1, 7)$, let $C=1-P-7-8-9-3-2-1$ and $C'=1-P-7-6-5-4-3-2-1$. If $(s, t)=(1, 6)$, let $C=1-P-6-5-4-3-2-1$ and $C'=1-P-6-7-8-9-3-2-1$. If $(s, t)=(1, 5)$, let $C=1-P-5-4-3-2-1$ and $C'=1-P-5-4-3-9-10-11-12-1$. The case $(s, t)=(1, 4)$ can be treated with a similar argument to $(s, t)=(1, 10)$. If $(s, t)=(1, 3)$, let  $C=1-P-3-9-10-11-12-1$ and $C'=1-P-3-9-8-7-6-12-1$. One of $C$ and $C'$ must be a big odd hole.

Suppose that $s=12$. The case where $t\in \{2, 4, 5, 7, 8, 10\}$ can be treated similarly as above. By symmetry it remains to deal with $(s, t)=(12, 9)$. now, we have $l(P)\ge4$, and so either $12-P-9-3-2-1-12$ or $12-P-9-3-4-5-6-12$ is a big odd hole.

Now, suppose that $(s, t)=(1, 11)$. Except 1,11 and possibly 12, no vertex of $H$ may have neighbors in $P^*$.  Since $N_H(F)\not\subseteq\{1, 11, 12\}$, there is a vertex, say $u$, in $V(H)\setminus\{1, 11, 12\}$ and an induced $(u, P^*)$-path $Q$ with interior in $V(F)$. We choose such a $Q$ that minimizes $\l(Q)$. By symmetry, we may assume that $u\in\{6, 7, 8, 9, 10\}$.

Let $R$ be an induced $1u$-path with $R^*\subseteq P^*\cup Q^*$. It is certain that no vertex of $H - \{1, 11 12, u\}$ may have neighbors in $R^*$. If $u=8$ then $l(R)\ge 4$ as otherwise $G[V(H)\cup R^*\setminus\{10, 11\}]={\cal P}$, let $C=1-R-8-9-3-2-1$ and $C'=1-R-8-7-6-5-4-3-2-1$. If $u=7$, let $C=1-R-7-8-9-3-2-1$ and $C'=1-R-7-6-5-4-3-2-1$. If $u=6$, let $C=1-R-6-5-4-3-2-1$ and $C'=1-R-6-7-8-9-3-2-1$. Each case leads to a contradiction as one of $C$ and $C'$ must be a big odd hole. Thus, we have that $9\le u\le 10$.

If $u=9$ and $N_{R^*}(12)\neq\emptyset$, let $R'$ be the shortest induced $(9, 12)$-path with interior in $R^*$, then either $12-R'-9-3-2-1-12$ or $12-R'-9-3-4-5-6-12$ is a big odd hole.
If $u=9$ and $N_{R^*}(12)=\emptyset$, then either $1-R-9-8-7-6-12-1$ or $1-R-9-3-4-5-6-12-1$ is a big odd hole.
If $u=10$ and $N_{R^*}(12)\neq\emptyset$, let $R'$ be the shortest induced $(10, 12)$-path with interior in $R^*$, then either $12-R'-10-9-3-2-1-12$ or $12-R'-10-9-3-4-5-6-12$ is a big odd hole. If $u=10$ and $N_{R^*}(12)=\emptyset$, then either $1-R-10-9-3-2-1$ or $1-R-10-9-3-4-5-6-12-1$ is a big odd hole. This proves Lemma~\ref{theo-1-5}.  \qed

\section{Jumps}

In this section, we always let $G$ be a heptagraph, and let $C=c_1\cdots c_7c_1$ be a 7-hole in $G$.  Since $G$ has no short cycles, we have that every vertex in $V(G)\setminus V(C)$ has at most one neighbor in $V(C)$. For two disjoint subsets $X$ and $Y$ of $V(G)$, we say that $X$ is {\em anticomplete} to $Y$ if no vertex of $X$ has neighbors in $Y$.

Let $sxyt$ be a segment of $C$. An induced $st$-path $P$ with $P^*\subseteq V(G)\setminus V(C)$ is called an {\em st e-jump across} $xy$. If $P$ is an $st$ $e$-jump across $xy$ such that $V(C)\setminus\{s, t, x, y\}$ is anticomplete to $P^*$, then we call $P$ a {\em local e-jump}. A local $e$-jump of length four is called a {\em short e-jump}.

Let $sct$ be a segment of $C$. An induced $st$-path $P$ with $P^*\subseteq V(G)\setminus V(C)$ is called an {\em st v-jump across} $c$. If $P$ is an $st$ $v$-jump across $c$ such that $V(C)\setminus\{c, s, t\}$ is anticomplete to $P^*$, then we call $P$ a {\em local v-jump}. A local $v$-jump of length five is called a {\em short v-jump}.

All $v$-jumps and $e$-jumps of $C$ are referred to as  {\em jumps} of $C$. A {\em non-local jump} is one which is not local. Clearly,

\renewcommand{\baselinestretch}{1}
\begin{itemize}
	
\item all $e$-jumps have length at least four, and all $v$-jumps have length at least five,

\item a jump $P$ is short if and only if  $V(C)\setminus V(P)$ is anticomplete to $P^*$, and

\item each local $e$-jump (resp. $v$-jump) has even (resp. odd) length.
\end{itemize}\renewcommand{\baselinestretch}{1.2}

\renewcommand{\baselinestretch}{1}	
\begin{lemma}\label{coro-2-1}
Suppose that, for $1\le i,j\le 7$, $C$ has a short $v$-jump across $c_i$ and a short $v$-jump across $c_j$ such that $d_C(c_i, c_j)\in \{2, 3\}$. Then $G$ induces a ${\cal P}$.
\end{lemma}\renewcommand{\baselinestretch}{1.2}	
\pf We first discuss the case where $d_C(c_i,c_j)=2$. Without loss of generality, suppose that $i=2$ and $j=7$. Let $P_1=c_1a_1a_2a_3a_4c_3$ and $P_2=c_1b_1b_2b_3b_4c_6$. To avoid a big odd hole on  $V(P_1)\cup V(P_2)\cup \{c_4,c_5\}$, there must exist $u\in P_1^*$ and $v\in P_2^*$ such that $u\sim v$.

If $u=a_4$, then $a_4b_4c_6c_5c_4c_3a_4$ is a 6-hole  when $v=b_4$,  $V(C)\cup\{b_1,b_2,b_3,b_4,a_4\}$ induces a ${\cal P}$ when $v=b_3$, $a_4b_2b_1c_1c_2c_3a_4$ is a 6-hole when $v=b_2$, and $a_4b_1c_1c_2c_3a_4$  is a 5-hole when $v=b_1$. Thus, we assume that $u\ne a_4$, and $v\ne b_4$ by symmetry.

If $u=a_3$, then $a_3b_3b_4c_6c_7c_1c_2c_3a_4a_3$ is a 9-hole  when $v=b_3$, $a_3b_2b_3b_4c_6c_5c_4c_3a_4a_3$ is a 9-hole when $v=b_2$, and $a_3b_1c_1c_2c_3a_4a_3$ is a 6-hole when  $v=b_1$. Thus, $u\ne a_3$, and $v\ne b_3$ by symmetry. Consequently, we have that $u\notin \{a_3,a_4\}$ and $v\notin \{b_3,b_4\}$.

If $u=a_2$, then $a_2b_2b_3b_4c_6c_7c_1c_2c_3a_4a_3a_2$ is an 11-hole  when $v=b_2$, and $a_2b_1b_2b_3b_4c_6c_5c_4c_3a_4a_3a_2$ is an 11-hole  when $v=b_1$.  If $(u,v)=(a_1,b_1)$, then $a_1=b_1$ to avoid a triangle $c_1a_1b_1c_1$, which implies an 11-hole  $a_1b_2b_3b_4c_6c_5c_4c_3a_4a_3a_2a_1$. Therefore,  the lemma holds if $d_C(c_i,c_j)=2$.

\medskip

Now, suppose $d_C(c_i,c_j)=3$. Without loss of generality, suppose that $i=3$ and $j=7$. Let $P_1=c_2a_1a_2a_3a_4c_4$ and $P_2=c_1b_1b_2b_3b_4c_6$. To avoid a 13-hole, there must exist $u\in P_1^*$ and $v\in P_2^*$ with $u\sim v$.

If $u=a_4$,  then $a_4b_4c_6c_5c_4a_4$ is a 5-hole  when $v=b_4$, $a_4b_3b_4c_6c_5c_4a_4$ is a 6-hole when $v=b_3$, $V(C)\cup\{b_4,b_3,b_2,b_1,a_4\}$ induces a ${\cal P}$  when $v=b_2$, and $a_4b_1c_1c_2c_3c_4a_4$ is a 6-hole  when $v=b_1$. Without loss of generality, we assume that $u\ne a_4$ and $v\ne b_4$.

If $u=a_3$, then $a_3b_3b_4c_6c_7c_1c_2a_1a_2a_3$ is a 9-hole  when $v=b_3$, $a_3b_2b_3b_4c_6c_7c_1c_2c_3c_4a_4a_3$ is an 11-hole when $v=b_2$, and $a_3b_1c_1c_2a_1a_2a_3$ is  a 6-hole when $v=b_1$. Now, suppose that $u\ne a_3$ and $v\ne b_3$, and this implies a short cycle in $G[c_1,c_2,a_1,a_2,b_1,b_2]$. Therefore, Lemma~\ref{coro-2-1} holds. \qed

\renewcommand{\baselinestretch}{1}	
\begin{lemma}\label{coro-2-2}
Let $P_1$ be a short $v$-jump, and $P_2$ be a short $e$-jump of $C$. If $P_1$ and $P_2$ share exactly one common end, and the other ends of them are not adjacent, then $G$ induces a ${\cal P}$.
\end{lemma}\renewcommand{\baselinestretch}{1.2}	
\pf Without loss of generality, suppose that $P_1=c_1a_1a_2a_3a_4c_3$ is a short $v$-jump across $c_2$, and $P_2=c_1b_1b_2b_3c_5$ is a short $e$-jump across $c_6c_7$. To avoid a n 11-hole on $V(P_1)\cup V(P_2)\cup \{c_4\}$, there must exist $u\in P_1^*$ and $v\in P_2^*$ with $u\sim v$.

To avoid a short cycle, $u$ cannot be $a_4$. If $u=a_3$, then $a_3b_3c_5c_4c_3a_4a_3$ is a 6-hole  when  $v=b_3$, $G[V(C)\cup\{a_4,a_3,a_1,b_2,b_3\}]={\cal P}$  when $v=b_2$, and $a_3b_1c_1c_2c_3a_4a_3$ is a 6-hole when $v=b_1$. If $u=a_2$, then $a_2b_2b_3c_5c_6c_7c_1c_2c_3a_4a_3a_2$ is an 11-hole when $v=b_2$, $a_2b_1b_2b_3c_5c_4c_3a_4a_3a_2$ is a 9-hole when $v=b_1$. If $u=a_1$, then $a_1b_3c_5c_6c_7c_1a_1$ is a 6-hole when $v=b_3$, $a_1b_2b_3c_5c_4c_3a_4a_3a_2a_1$ is a 9-hole  when $v=b_2$. If $(u,v)=(a_1,b_1)$, then $a_1=b_1$ to avoid a triangle $c_1a_1b_1c_1$, and so $a_1b_2b_3c_5c_4c_3a_4a_3a_2a_1$ is a 9-hole. Therefore, Lemma~\ref{coro-2-2} holds. \qed

\renewcommand{\baselinestretch}{1}	
\begin{lemma}\label{coro-2-3}
If $C$ has two short $e$-jumps sharing exactly one common end, then $G$ induces a ${\cal P}'$.
\end{lemma}\renewcommand{\baselinestretch}{1.2}	
\pf Without loss of generality, suppose that $P_1=c_1a_1a_2a_3c_4$ and $P_2=c_1b_1b_2b_3c_5$ are two short $e$-jumps. To avoid a 9-hole on $V(P_1)\cup V(P_2)$, there must exist $u\in P_1^*$ and $v\in P_2^*$ with $u\sim v$.

If $u=a_3$ or $v=b_3$, then a short cycle occurs. Thus, we have that $u\ne a_3$ and $v\ne b_3$. If $u=a_2$, then we have a 6-hole $a_2b_2b_3c_5c_4a_3a_2$ when $v=b_2$, and an induced ${\cal P}'$ on $V(C)\cup\{b_1,b_2,b_3,a_2,a_3\}$ when $v=b_1$. The same contradiction occurs if $v=b_2$. If $(u,v)=(a_1,b_1)$ then $a_1=b_1$ to avoid a triangle $c_1a_1b_1c_1$, and so $G[V(C)\cup\{a_1,a_2,a_3,b_2,b_3\}]={\cal P}'$. This proves Lemma~\ref{coro-2-3}. \qed

\renewcommand{\baselinestretch}{1}	
We say that
\begin{itemize}
\item $C$ is of {\em type $1$} if it has two local $v$-jumps sharing  exactly one common end,

\item $C$ is of {\em type $2$} if $C$ is not of {\em type $1$}, and has a local $v$-jump $P_1$ and a local $e$-jump $P_2$ such that $P_1$ and $P_2$ share exactly one common end and the other ends of them are not adjacent, and

\item $C$ is of {\em type $3$} if $C$ is not of {\em type $1$} or {\em type $2$}, and has two local $e$-jumps sharing exactly one common end.
\end{itemize}\renewcommand{\baselinestretch}{1.2}	

The following summations of subindexes are taken modulo 7, and we set $7+1\equiv 1$.

\renewcommand{\baselinestretch}{1}
\begin{lemma}\label{lem-3-1}
Suppose that $G$ induces no ${\cal P}$ and ${\cal P}'$, and suppose that $C$ is of type $1$ with two local $v$-jumps $P_1$ and $P_2$ that share exactly one common end $c_j$ for some $j\in\{1, 2, \ldots, 7\}$. Then, at least one of $P_1$ and $P_2$ is not short, and $C$ has a short jump $T$ with interior in $P_1^*\cup P_2^*$ such that
\begin{itemize}
\item [$(a)$]  $T$ is a $v$-jump across $c_j$, or an $e$-jump across $c_{j-1}c_j$ or $c_jc_{j+1}$, and

\item [$(b)$] none of $P_1$ and $P_2$ is short if $T$ is a short $v$-jump across $c_j$.
\end{itemize}
\end{lemma}\renewcommand{\baselinestretch}{1.2}
\pf It follows directly from Lemmas~\ref{coro-2-1} that at least one of $P_1$ and $P_2$ is not short. The statement $(b)$ follows directly from the fact that $G$ has no short cycles. Now it is left to prove $(a)$.

Without loss of generality, suppose that $j=1$, and suppose that for each short jump $Q$ of $C$ with interior in $P_1^*\cup P_2^*$, $Q$ is neither a $v$-jump across $c_1$ nor an $e$-jump across $c_1c_{2}$ or $c_{1}c_7$. We may choose $P_1$ and $P_2$ such that $|P_1^*\cup P_2^*|$ is minimum. Let $P_1=c_1a_1a_2a_3\ldots a_kc_3$ and $P_2=c_1b_1b_2b_3\ldots b_tc_6$. Let $D_1=V(P_1[a_4,a_{k-1}])$ and $D_2=V(P_2[b_4,b_{t-1}])$.

\begin{claim}\label{clm-3-1}
$D_1\cup \{a_k\}$ is  disjoint from and anticomplete to $D_2\cup \{b_t\}$.
\end{claim}	
\pf Since $P_1$ and $P_2$ are both local, we have that $a_k\notin V(P_2)$ and $b_t\not\in V(P_1)$, and so $a_k\not\sim b_t$ to avoid a short cycle. Suppose that the claim is not true, and suppose by symmetry that there is an $(a_k, D_2)$-path in $G[D_1\cup D_2\cup\{a_k\}]$.
Thus $P_2$ cannot be short, and so $N_{D_2}(c_7)\ne\emptyset$. We may choose $P'$ to be a $(c_7, \{c_2,c_3\})$-path   with shortest length and interior in $D_1\cup D_2\cup \{a_k\}$. Let $x$ be the end of $P'$ other than $c_7$. It is certain that $V(C)\setminus\{c_7, x\}$ is anticomplete to $P'^*$, which implies that $P'$ is either a short $v$-jump across $c_1$ or a short $e$-jump across $c_1c_2$, contradicting our assumption. This proves Claim ~\ref{clm-3-1}.  \qed

Note that both $\l(P_1)$ and $\l(P_2)$ are odd. If $a_3=b_3$, then $G[D_1\cup D_2\cup \{a_k,b_t,a_3,c_3,c_4,c_5,c_6\})]$ is a 7-hole, and so $P_1$ and $P_2$ are both short, contradicting Lemma~\ref{coro-2-1}. Hence, $a_3\ne b_3$. To avoid short cycles, we have $a_3\not\sim b_3$. To avoid big odd holes, we have $\{a_1, a_2\}\cap \{b_1, b_2\}=\emptyset$. By the minimality of $|P_1^*\cup P_2^*|$, we have that $P_1^*$ is disjoint from and anticomplete to $P_2^*$. This implies that $G[P_1\cup P_2\cup\{c_4,c_5\}]$ is a big odd hole. This proves Lemma~\ref{lem-3-1}. \qed

\renewcommand{\baselinestretch}{1}
\begin{lemma}\label{lem-3-2}
Suppose that $G$ induces no ${\cal P}$ and ${\cal P}'$, and suppose that $C$ is of type $2$ with a local $v$-jump $P_1$ and a local $e$-jump $P_2$ such that $P_1$ and $P_2$ share exactly one  common end $c_j$ for some $j\in\{1, 2, \ldots, 7\}$. Then, $C$ has a short jump $T$ with interior in $P_1^*\cup P_2^*$ such that
\begin{itemize}
\item [$(a)$]  $T$ is either a $v$-jump across $c_j$, or an $e$-jump across $c_{j-1}c_j$ or $c_jc_{j+1}$, and

\item [$(b)$] $P_2$ is not short, and $P_1$ is not short if $T$ is a short $v$-jump across $c_j$.
\end{itemize}
\end{lemma}\renewcommand{\baselinestretch}{1.2}
\pf We only need to prove $(a)$. The statement $(b)$ follows from $(a)$ directly.

Suppose that $(a)$ does not hold. Without loss of generality, we suppose that $j=1$, and suppose that $P_1$ and $P_2$ are chosen with $|P_1^*\cup P_2^*|$ minimum. Let $P_1=c_1a_1a_2\dots a_kc_3$ be a local $v$-jump, and $P_2=c_1b_1b_2\dots b_tc_5$ be a local  $e$-jump.  Let $D_1=V(P_1[a_3,a_k])$ and $D_2=V(P_2[b_3,b_t])$. Since $C$ is not of type 1, we have that $N_{P_2}(c_6)=\emptyset$. Furthermore, we have the following
\begin{claim}\label{clm-3-2}
	$D_1\cup \{a_k\}$ is disjoint from and anticomplete to $D_2\cup \{b_t\}$.
\end{claim}	
\pf Since both $P_1$ and $P_2$ are local, we have that $a_k\ne b_t$, and so $a_k\not\sim b_t$ to avoid short cycles. If the claim is not true, then  $G[D_1\cup D_2\cup \{a_k\}]$ has an $(a_k, D_2)$-path, and so $C$ has a $v$-jump $P'$ across $c_4$. Since $P'$ is not local and $N_{P^*_2}(c_6)=\emptyset$, we have that $N_{D_2}(c_7)\ne\emptyset$, and we may choose $Q$ to be a  $(c_7, \{c_2,c_3\})$-path with shortest length and interior in $D_1\cup D_2\cup \{a_k\}$. Let $x$ be the end of $Q$ other than $c_7$. It is certain that $V(C)\setminus\{c_7, x\}$ is anticomplete to $Q^*$, which implies that $Q$ is either a short $v$-jump across $c_1$ or a short $e$-jump across $c_1c_2$, a contradiction. This proves Claim~\ref{clm-3-2}.

With the same arguments as that used in the proof of Lemma~\ref{lem-3-1}, we have that $P_1^*$ is disjoint from and anticomplete to $P_2^*$, which implies a big odd hole on $P_1\cup P_2\cup\{c_4\}$. This proves Lemma~\ref{lem-3-2}.\qed

\renewcommand{\baselinestretch}{1}
\begin{lemma}\label{lem-3-3}
Suppose that $G$ induces no ${\cal P}$ and ${\cal P}'$, and suppose that $C$ is of type $3$ with two local $e$-jumps $P_1$ and $P_2$ that share exactly one  common end $c_j$ for some $j\in\{1, 2, \ldots, 7\}$. Then
\begin{itemize}
\item [$(a)$] $C$ has a short jump with interior in $P_1^*\cup P_2^*$ which is either a $v$-jump across $c_j$, or an $e$-jump across $c_{j-1}c_j$ or $c_jc_{j+1}$, and

\item [$(b)$] none of $P_1$ and $P_2$ is short.
\end{itemize}
\end{lemma}\renewcommand{\baselinestretch}{1.2}
\pf We only need to prove $(a)$. Suppose to its contrary that $(a)$ is not true. Without loss of generality, suppose that $j=1$ and $P_1$, $P_2$ are chosen such that $|P_1^*\cup P_2^*|$ is minimum. Let $P_1=c_1a_1\dots a_kc_4$ and $P_2=c_1b_1\dots b_tc_5$. Let $D_1=V(P_1[a_2,a_k])$ and $D_2=V(P_2[b_2,b_t])$.

Since $C$ is not type 1 or type 2, we have that $N_{P^*_1}(c_3)=N_{P^*_2}(c_6)=\emptyset$. By Lemma~\ref{coro-2-3}, we have that one of $P_1$ and $P_2$, say $P_2$, is not short. Thus, $N_{P^*_2}(c_7)\neq \emptyset$.

\begin{claim}\label{clm-3-3}
$D_1\cup \{a_k\}$ is disjoint from and anticomplete to $D_2\cup \{b_t\}$.
\end{claim}		
\pf Since both $P_1$ and $P_2$ are local jumps, we have that $a_k\not\in P^*_2$ and $b_t\not\in P^*_1$, and so $a_k\not\sim b_t$ to avoid a 4-cycle $c_4a_kb_tc_5c_4$.

Suppose that $N_{P^*_2}(a_k)\ne\emptyset$. Choose $x\in N_{P^*_2}(a_k)$ to be a vertex closest to $b_t$. By the minimality of $|P_1^*\cup P_2^*|$, we may assume that $x=a_{k-1}$ and $P_1[c_1, x]=P_2[c_1, x]$. Notice that $N_{P^*_2}(c_7)\ne \emptyset$. If $N_{P^*_1}(c_2)\ne \emptyset$, we may choose $P$ to be an induced $c_2c_7$-path with $P^*\subseteq P_1^*\cup P_2^*\setminus\{a_1, a_k, b_t\}$, then $P^*$ is anticomplete to $V(C)\setminus\{c_2, c_7\}$, which gives a short jump across $c_1$, a contradiction. Thus, $N_{P^*_1}(c_2)=\emptyset$, which implies that $P_1$ is a short jump. Hence, $k=3$ and $c_7$ must have neighbors in $P_2[x, b_t]-x$. Consequently, $\l(P_2[x, b_t])$ is odd and at least 7. Thus $c_4c_5b_tP_2[x, b_t]xa_3c_4$ is a big odd hole.  Therefore, $N_{P^*_2}(a_k)=\emptyset$ and $N_{P^*_1}(b_t)=\emptyset$ by symmetry.

If the claim is not true, let $i$ and $j$ be the largest indexes such that $a_i\sim b_j$, then either $C$ has a short jump across $c_1$ when $P_1$ is not short, or  $c_4P_1[c_4, a_i]a_ib_jP_2[b_i, c_5]c_5c_4$ is a big odd hole when $P_1$ is short.  This proves Claim~\ref{clm-3-3}. \qed

With the similar arguments  as that used in  the proof of Lemma~\ref{lem-3-1}, we conclude that $P_1^*$ is disjoint from and anticomplete to $P_2^*$, which gives a big odd hole $G[V(P_1)\cup V(P_2)]$. Therefore, Lemma~\ref{lem-3-3} holds.\qed

\renewcommand{\baselinestretch}{1}	
\begin{lemma}\label{lem-3-4}
Let $P$ be a jump of $C$.  Suppose that $G$ induces no ${\cal P}$ and ${\cal P}'$ and $P$ is not a local jump. Then $C$ has a short jump with interior in $P^*$.
\end{lemma}\renewcommand{\baselinestretch}{1.2}	
\pf Without loss of generality, suppose that  $P$ is a $v$-jump across $c_2$ or an $e$-jump across $c_2c_3$.

If $P$ is a $v$-jump and $N_{P^*}(c_2)\ne\emptyset$, then let $Q$ be a  ($c_2, \{c_4, c_5, c_6, c_7\}$)-path with shortest length and $Q^*\subseteq P^*$. If $P$ is an $e$-jump and $N_{P^*}(c_2)\cup N_{P^*}(c_{3})\ne\emptyset$, then let $Q$ be a  ($\{c_2, c_3\}, \{c_5, c_6, c_7\}$)-path with shortest length and $Q^*\subseteq P^*$. It is easy to verify that $Q$ must be a short jump in both cases.

Now suppose that $N_{P^*}(c_2)=\emptyset$ when $P$ is a $v$-jump, and $N_{P^*}(c_2)\cup N_{P^*}(c_{3})=\emptyset$ when $P$ is an $e$-jump.
We only prove the case where $P$ is a $v$-jump. The case that $P$ is an $e$-jump can be treated with almost the same arguments.

Suppose to the contrary that the lemma is not true. Firstly, we show that
\begin{equation}\label{eqa-nonlocal-c4-0}
N_{P^*}(c_4)=\emptyset.
\end{equation}

Suppose that $N_{P^*}(c_4)\ne\emptyset$. Let $Q$ be the shortest $(c_1, c_4)$-path with $Q^*\subseteq P^*$. Then $N_{Q^*}(c_2)=N_{Q^*}(c_{3})=\emptyset$ and $N_{Q^*}(\{c_5, c_6, c_7\})\ne\emptyset$.

If $N_{Q^*}(c_6)\ne \emptyset$, then let $Q_{1, 6}$ be the shortest $c_1c_6$-path and $Q_{4, 6}$ be the shortest $c_4c_6$-path, both with interior in $Q^*$. If $Q_{1, 6}$ and $Q_{4, 6}$ are both local, then by applying Lemma~\ref{lem-3-1} to $Q_{1, 6}$ and $Q_{4, 6}$, we can find a short jump as required. Thus by symmetry we assume that $Q_{1, 6}$ is not local. Then, $N_{Q^*_{1, 6}}(c_5)\ne\emptyset$. Thus, either $C$ has a short $e$-jump across $c_6c_7$ when $N_{Q^*_{1, 6}}(c_7)=\emptyset$, or $C$ has a short $v$-jump  across $c_6$ when $N_{Q^*_{1, 6}}(c_7)\ne \emptyset$. This shows that $N_{Q^*}(c_6)=\emptyset$.

If $N_{Q^*}(c_5)\ne \emptyset$ and $N_{Q^*}(c_7)\ne \emptyset$, then the shortest $c_5c_7$-path, with interior in $Q^*$, is a short $v$-jump as required. Otherwise, we may assume by symmetry that $N_{Q^*}(c_5)\ne \emptyset$ and $N_{Q^*}(c_7)=\emptyset$, then the shortest $c_1c_5$-path, with interior in $Q^*$, is a short $e$-jump as required. Therefore, (\ref{eqa-nonlocal-c4-0}) holds.

By symmetry, we may suppose that  $N_{P^*}(c_4)=\emptyset$ and $N_{P^*}(c_7)=\emptyset$. Thus, a $(\{c_1, c_3\}, \{c_5, c_6\})$-path, with shortest length and interior in $Q^*$, is a short jump as required. This proves  Lemma~\ref{lem-3-4}.  \qed


\section{Proof of Theorem~\ref{theo-1-3}}

In this section we prove Theorem~\ref{theo-1-3}. If a heptagraph has no 7-hole, then it is bipartite. Thus we always use $G$ to denote a heptagraph, use $C=c_1\cdots c_7c_1$ to denote a 7-hole in $G$, and let ${\cal X}$ be the set of all vertices which are in the interior of some short jumps of $C$. For two integers $i$ and $j$ with $1\le i<j\le 7$, we use $X_{i, j}$ to denote the set of all vertices which are in the interior of some short jumps joining $c_i$ and $c_j$.

The proof of Theorem~\ref{theo-1-3} is divided into a several  lemmas. By Lemma~\ref{lem-2-1}, we have that $C$ must have some local jumps. We say that two local jumps are {\em equivalent} if they have the same ends. We start from the case that all local jumps of $C$ are equivalent. After that, we discuss the cases where $C$ is of type $i$ for some $i\in\{1, 2, 3\}$. At last we consider the case where $C$ has two kinds of equivalent local jumps and is not of  type $i$ for any $i$.

In the proof of each lemma, we will choose a subset ${\cal D}$ of $V(G)$ which is disjoint from $V(C)\cup {\cal X}$, and call a short jump  {\em bad} if it has some interior vertex in ${\cal D}$. It is certain that
\begin{equation}\label{eqa-3-3-1}
\mbox{$C$ has no bad jumps.}
\end{equation}

We always use ${\cal N}$ to denote the set of vertices in $V(C)\cup {\cal X}$ that have neighbors in ${\cal D}$.

In the proofs, Lemmas~\ref{lem-3-1}, \ref{lem-3-2}, \ref{lem-3-3} and \ref{lem-3-4} will be cited frequently. We use Lemma~\ref{lem-3-4}($P$) to denote the set of short jumps obtained by applying Lemma~\ref{lem-3-4} to a jump $P$, and use Lemmas~\ref{lem-3-1}$(P, Q)$ to denote the set of short jumps obtained by applying Lemmas~\ref{lem-3-1} to local $v$-jumps $P$ and $Q$ which share exactly one end. Similarly, we define Lemma~\ref{lem-3-2}$(P, Q)$, and Lemma~\ref{lem-3-3}$(P, Q)$.

Since $G$ has no triangles, we have that each clique cutset is a single vertex or the two ends of an edge, which is a parity-star cutset. In the rest of the paper, we always choose $G$ to be a heptagraph such that
\renewcommand{\baselinestretch}{1}	
\begin{itemize}
\item $\delta(G)\ge 3$, $G$ induces no ${\cal P}$ and ${\cal P}'$, and $G$ has no clique cutsets and no $P_3$-cutsets.
\end{itemize}\renewcommand{\baselinestretch}{1.2}	

\begin{lemma}\label{lem-unique-local}
Suppose that all local jumps of $C$ are equivalent. Then $G$ admits a parity star-cutset.
\end{lemma}	
\pf  Let $P$ be a local jump of the shortest length.
We first suppose that $P$ is a local $v$-jump across $c_2$.  Then ${\cal X}=X_{1, 3}$.  Since $N_{\cal X}(c_7)=\emptyset$ and $d(c_7)\ge 3$, we choose ${\cal D}$ to be the vertex set of a maximal connected induced subgraph such that $N_{{\cal D}}(c_7)\ne \emptyset$ and ${\cal D}\cap (V(C)\cup {\cal X})=\emptyset$. Then (\ref{eqa-3-3-1}) holds.

Let $X_1=N(c_1)\cap X_{1, 3}$ and $X_3=X_{1, 3}\setminus X_1$.
	
Suppose that $(X_3\cup\{c_3\})\cap {\cal N}\ne\emptyset$. Let $Q_{3, 7}$ be a  $c_3c_7$-path with shortest length and $Q^*_{3, 7}\subseteq {\cal D}\cup X_3$. Since $Q_{3, 7}$ is not a local jump, we have that $N_{Q^*_{3, 7}}(\{c_4, c_5, c_6\})\ne \emptyset$, and so Lemma~\ref{lem-3-4}($Q_{3, 7}$) has a bad jump, contradicting (\ref{eqa-3-3-1}). Therefore, $(X_3\cup\{c_3\})\cap {\cal N}=\emptyset$.

Suppose that $c_4\in {\cal N}$. Let $Q_{4, 7}$ be a $c_4c_7$-path with shortest length and interior in ${\cal D}$. Since $Q_{4, 7}$ is not a local jump, we have that Lemma~\ref{lem-3-4}($Q_{4, 7}$) has a bad jump, contradicting (\ref{eqa-3-3-1}). Therefore, $c_4\not\in {\cal N}$.

With a similar argument we can show that $c_5\not\in {\cal N}$ and $c_6\not\in {\cal N}$. Thus ${\cal N}\subseteq X_1\cup\{c_1, c_2, c_7\}$.

Since $G$ has no $P_3$-cutsets, we have that ${\cal N}\cap X_1\ne\emptyset$. If $c_2\in {\cal N}$, then the shortest $c_2c_7$-path, with interior in ${\cal D}\cup\{x_1\}$, is a local jump, a contradiction. Therefore, ${\cal N}\subseteq X_1\cup \{c_1,c_7\}$. Since every two vertices in $X_1\cup\{c_7\}$ are joined by an induced path of length six or eight with interior in $X_1\cup\{c_3, c_4,c_5,c_6\}$, we have that $X_1\cup\{c_1,c_7\}$ is a parity star-cutset.

\medskip

Next, we suppose that $P$ is a local $e$-jump  across $c_2c_3$. Then, ${\cal X}=X_{1, 4}$. Let $X_1=N(c_1)\cap X_{1, 4}$ and $X_4=X_{1, 4}\setminus X_1$.

Since all local jumps of $C$ are equivalent to $P$, we have that $N_{{\cal X}}(\{c_2, c_3, c_5, c_6, c_7\})=\emptyset$. Since  $d(c_6)\ge 3$, we choose ${\cal D}$ to be the vertex set of a maximal connected induced subgraph such that $N_{{\cal D}}(c_6)\ne \emptyset$ and ${\cal D}\cap (V(C)\cup {\cal X})=\emptyset$. Thus (\ref{eqa-3-3-1}) still holds. We claim that
\begin{equation}\label{eqa-c1c2c3c4-N}
(X_4\cup\{c_4\})\cap {\cal N}=(X_1\cup\{c_1\})\cap {\cal N}=\{c_2, c_3\}\cap {\cal N}=\emptyset.
\end{equation}

Suppose that $(X_4\cup\{c_4\})\cap {\cal N}\ne\emptyset$. Let $Q_{4, 6}$ be a  $c_4c_6$-path with shortest length and $Q^*_{4, 6}\subseteq {\cal D}\cup X_4$. Since $Q_{4, 6}$ is not a local jump, we have that $N_{Q^*_{4, 6}}(\{c_1, c_2, c_3,c_7\})\ne\emptyset$, and so Lemma~\ref{lem-3-4}($Q_{4, 6}$) has a bad jump, contradicting (\ref{eqa-3-3-1}). Hence, $(X_4\cup\{c_4\})\cap {\cal N}=\emptyset$, and $(X_1\cup\{c_1\})\cap {\cal N}=\emptyset$ by symmetry.

If $c_3\in {\cal N}$, let $Q_{3, 6}$ be a $c_3c_6$-path with shortest length and $Q^*_{3, 6}\subseteq {\cal D}$, then $Q_{3, 6}$ is not a local jump, and $N_{Q^*_{3, 6}}(\{c_1, c_2, c_7\})\ne\emptyset$. Consequently, $C$ has a bad jump in Lemma~\ref{lem-3-4}($Q_{3, 6}$). Therefore, $c_3\not\in {\cal N}$, and $c_2\not\in {\cal N}$ by symmetry.
This proves (\ref{eqa-c1c2c3c4-N}).

By (\ref{eqa-c1c2c3c4-N}), ${\cal N}\subseteq \{c_5,c_6,c_7\}$ and induces a $P_3$-cutset of $G$, contradicting the choice of $G$. This completes the proof of Lemma~\ref{lem-unique-local}. \qed

\medskip

Now suppose that $C$ has at least two kinds of equivalent local jumps. If $C$ is of type $i$ for some $i\in \{1, 2, 3\}$, then we always choose $j$ and the two local jump $P_1$ and $P_2$ such that $P_1$ and $P_2$ share $c_j$ and $|P_1^*\cup P_2^*|$ is minimum. Without loss of generality, suppose that $j=1$.


\begin{lemma}\label{theo-3-1}
Suppose that $C$ is of type $1$. Then $G$ admits a parity star-cutset.
\end{lemma}	
\pf Since $P_1$ and $P_2$ are local jumps sharing $c_1$, we have that Lemma~\ref{lem-3-1}($P_1, P_2$) has a short jump $T$, with $T^*\subseteq P_1^*\cup P_2^*$, which is a $v$-jump across $c_1$, or an $e$-jump across $c_1c_2$ or  $c_1c_7$.

Firstly, we prove
\begin{claim}\label{clm-3-4}
Lemma~{\em \ref{theo-3-1}} holds if $T$  is a short $v$-jump across $c_1$.
\end{claim}		
\pf Suppose that $T$ is a $v$-jump across $c_1$. It is certain that both $P_1$ and $P_2$ are not short.

By the minimality of $|P_1^*\cup P_2^*|$, $C$ has no short $v$-jumps across $c_2$ or $c_7$. Since the two ends of any jump of $C$ are not adjacent, by Lemmas~\ref{coro-2-1}, \ref{coro-2-2}, \ref{lem-3-1}, and \ref{lem-3-2}, we have that no short jumps have end $c_4$ or $c_5$. Thus, ${\cal X}=X_{2,7}\cup X_{2,6}\cup X_{3,7}$, and $c_4$ has no neighbor in ${\cal X}\cup \{c_1,c_2,c_6,c_7\}$ as each vertex in $N_{{\cal X}}(c_4)$ provides us with a short jump starting from $c_4$.

Since $d(c_4)\ge 3$, we choose ${\cal D}$ to be the vertex set of a maximal connected induced subgraph such that $N_{{\cal D}}(c_4)\ne \emptyset$ and ${\cal D}\cap (V(C)\cup {\cal X})=\emptyset$. It is certain that (\ref{eqa-3-3-1}) holds.

Let $X_3=X_{3,7}\cap N(c_3)$. First we claim that
\begin{equation}\label{eqa-3-4}
\mbox{$((X_{2,6}\cap N(c_2))\cup (X_{2,7}\cap N(c_2))\cup\{c_2\})\cap {\cal N}=\emptyset$.}
\end{equation}	

Suppose that (\ref{eqa-3-4}) is not true. Let $Q_{2, 4}$ be a $c_2c_4$-path with shortest length and $Q^*_{2, 4}\subseteq {\cal D}\cup (X_{2,6}\cap N(c_2))\cup (X_{2,7}\cap N(c_2))$. If $Q_{2, 4}$ is a local jump then Lemma~\ref{lem-3-1}($Q_{2, 4}, T)$ has a bad jump. If $Q_{2, 4}$ is not a local jump then Lemma~\ref{lem-3-4}($Q_{2, 4}$) has a bad jump. Both contradict (\ref{eqa-3-3-1}). Therefore, (\ref{eqa-3-4}) holds.

Next we claim that
\begin{equation}\label{eqa-3-5}
\mbox{$c_1\notin {\cal N}$.}
\end{equation}		

Suppose that $c_1\in {\cal N}$. Let $Q_{1, 4}$ be a $c_1c_4$-path with shortest  length and interior in ${\cal D}$. By (\ref{eqa-3-4}), we have that $N_{Q^*_{1, 4}}(c_2)=\emptyset$.
If $Q_{1, 4}$ is a local $e$-jump, then Lemma~\ref{lem-3-2}($P_2, Q_{1, 4}$) has a bad jump. If $Q_{1, 4}$ is not local, then Lemma~\ref{lem-3-4}($Q_{1, 4}$) has a bad jump. Both contradict (\ref{eqa-3-3-1}). Therefore, (\ref{eqa-3-5}) holds.

Now we claim that
\begin{equation}\label{eqa-3-6}
\mbox{$((X_{2,7}\setminus N(c_2))\cup (X_{3,7}\setminus N(c_3))\cup\{c_7\})\cap {\cal N}=\emptyset$.}
\end{equation}		

Suppose that (\ref{eqa-3-6}) does not hold. Let $Q_{4, 7}$ be a $c_4c_7$-path with shortest length and interior in ${\cal D}\cup (X_{2,7}\setminus N(c_2))\cup (X_{3,7}\setminus N(c_3))$. By (\ref{eqa-3-4}) and (\ref{eqa-3-5}), we have that $N_{Q^*_{4, 7}}(\{c_1, c_2\})=\emptyset$. If $Q_{4, 7}$ is a local jump, then Lemma~\ref{lem-3-2}($Q_{4, 7}, T$) has a bad $e$-jump across $c_1c_7$. If $Q_{4, 7}$ is not local, then $N_{Q^*_{4, 7}}(c_3)\ne\emptyset$, and the shortest ($c_3, \{c_5, c_6, c_7\})$-path, with interior in $Q^*_{4, 7}$, is a bad jump. Both contradict (\ref{eqa-3-3-1}). Therefore, (\ref{eqa-3-6}) holds.

Finally we claim that
\begin{equation}\label{eqa-3-7}
\mbox{$((X_{2,6}\setminus N(c_2))\cup\{c_6\})\cap {\cal N}=\emptyset$.}
\end{equation}	

Suppose it is not true.  Let $Q_{4, 6}$ be a $c_4c_6$-path with shortest length and interior in ${\cal D}\cup (X_{2,6}\setminus N(c_2))$. By (\ref{eqa-3-4}), (\ref{eqa-3-5}) and (\ref{eqa-3-6}), we have that $N_{Q^*_{4, 6}}(\{c_1, c_2, c_7\})=\emptyset$. If $N_{Q^*_{4, 6}}(c_3)=\emptyset$ then $Q^*_{4, 6}$ is a local $v$-jump, and Lemma~\ref{lem-3-1}($P_2, Q_{4, 6}$) has a bad jump (with the end either $c_4$ or $c_5$). Otherwise, the shortest ($c_3, \{c_5, c_6\})$-path, with interior in $Q^*_{4, 6}\subseteq {\cal D}\cup (X_{2,6}\setminus N(c_2))$, is a bad jump. Both contradict (\ref{eqa-3-3-1}). Therefore, (\ref{eqa-3-7}) holds.

\medskip

By  (\ref{eqa-3-4}),  (\ref{eqa-3-5}),  (\ref{eqa-3-6}), and (\ref{eqa-3-7}), we have that ${\cal N}\subseteq X_3\cup \{c_3,c_4,c_5\}$. Since $G$ has no $P_3$-cutsets, we have ${\cal N}\cap X_3\ne\emptyset$. If $N_{{\cal D}}(c_5)\ne \emptyset$, then there is a bad $v$-jump across $c_4$, contradicting (\ref{eqa-3-3-1}). Hence, we have that ${\cal N}\subseteq X_3\cup \{c_3,c_4\}$. Notice that each pair of vertices in $X_3\cup\{c_4\}$ are joined by an induced path of length 6 with interior in $X_{3,7}\setminus N(c_3)\cup\{c_5,c_6,c_7\}$. Thus  ${\cal N}\cup \{c_3,c_4\}$ is a parity star-cutset. This proves Claim~\ref{clm-3-4}.\qed	

\medskip

Now suppose that $C$ has no short $v$-jumps across $c_1$ with interior in $P_1^*\cup P_2^*$. Thus $T$ must be a short $e$-jump across either $c_1c_2$ or $c_1c_7$. Without loss of generality, suppose that $T$ is a short $e$-jump across $c_1c_2$.

\begin{claim}\label{clm-3-5}
Lemma~$\ref{theo-3-1}$ holds if $C$ has no short $v$-jump across $c_3$.
\end{claim}	
\pf Suppose that $C$ has no short $v$-jump across $c_3$. By Lemma~\ref{coro-2-1}, we have that
either $X_{1,3}=\emptyset$ or $X_{1,6}=\emptyset$.
Since the two ends of any jump of $C$ are not adjacent, by Lemmas~\ref{coro-2-1}, \ref{coro-2-2}, \ref{lem-3-1} and \ref{lem-3-2}, we have that no short jumps contain $c_4$ or $c_5$. Thus, ${\cal X}=X_{2,7}\cup X_{2,6}\cup X_{3,7}\cup X_{1,3}\cup X_{1,6}$, and $c_5$ has no neighbors in ${\cal X}\cup \{c_1,c_2,c_3,c_7\}$ as each vertex in $N_{{\cal X}}(c_5)$ provides us with a short jump starting from $c_5$.

Since $d(c_5)\ge 3$, we choose ${\cal D}$ to be the vertex set of a maximal connected induced subgraph such that $N_{{\cal D}}(c_5)\ne \emptyset$ and ${\cal D}\cap (V(C)\cup {\cal X})=\emptyset$. Then (\ref{eqa-3-3-1}) holds.

Let $X_6=N(c_6)\cap X_{2,6}$ and $X'_6=N(c_6)\cap X_{1,6}$. With the similar arguments as that used in the proof of Claim~\ref{clm-3-4}, we now prove that
\begin{equation}\label{eqa-3-9-0}
{\cal N}\subseteq (X_6'\cup X_6)\cup \{c_4,c_5,c_6\}.
\end{equation}

Suppose that $[(X_{3,7}\setminus N(c_7))\cup (X_{1,3}\cap N(c_3))\cup\{c_3\}]\cap {\cal N}\ne \emptyset$. Let $Q_{3, 5}$ be a $c_3c_5$-path with shortest length and $Q^*_{3, 5}\subseteq {\cal D}\cup (X_{3,7}\setminus N(c_7))\cup (X_{1,3}\cap N(c_3))$. If $Q_{3, 5}$ is a local jump, then  Lemma~\ref{lem-3-2}($Q_{3, 5}, T$) has a bad jump with end either $c_4$ or $c_5$. If $Q_{3, 5}$ is not a local jump, then Lemma~\ref{lem-3-4}($Q_{3, 5}$) has a bad jump. Both contradict (\ref{eqa-3-3-1}). Therefore,  ${\cal N}\subseteq (V(C)\setminus \{c_3\})\cup X_{2,7}\cup X_{2,6}\cup X_{1,6}\cup (X_{3,7}\cap N(c_7))\cup (X_{1,3}\setminus N(c_3))$.

Suppose that $[(X_{2,7}\setminus N(c_2))\cup (X_{3,7}\cap N(c_7))\cup\{c_7\}]\cap {\cal N}\ne \emptyset$.  Let $Q_{5, 7}$ be a $c_5c_7$-path with shortest length and $Q^*_{5, 7}\subseteq {\cal D}\cup (X_{2,7}\setminus N(c_2))\cup (X_{3,7}\cap N(c_7))$. Then, we have $N_{Q^*_{5, 7}}(c_3)=\emptyset$ as $c_3\not\in {\cal N}$. If $Q_{5, 7}$ is a local jump then Lemma~\ref{lem-3-2}($Q_{5, 7}, T$) has a bad jump. If $Q_{5, 7}$ is not a local jump then Lemma~\ref{lem-3-4}($Q_{5, 7}$) has a bad jump. Both contradict (\ref{eqa-3-3-1}). Thus we have ${\cal N}\subseteq (V(C)\setminus \{c_3, c_7\})\cup X_{2,6}\cup X_{1,6}\cup (X_{1,3}\setminus N(c_3))\cup (X_{2,7}\cap N(c_2))$.

Suppose that $[(X_{2,7}\cap N(c_2))\cup (X_{2,6}\setminus N(c_6))\cup\{c_2\}]\cap {\cal N}\ne \emptyset$. Let $Q_{2, 5}$ be a $c_2c_5$-path with shortest length and $Q^*_{2, 5}\subseteq {\cal D}\cup (X_{2,7}\cap N(c_2))\cup (X_{2,6}\setminus N(c_6))$.  Then  $N_{Q^*_{2, 5}}(c_3)=N_{Q^*_{2, 5}}(c_7)=\emptyset$ as $c_3, c_7\not\in {\cal N}$. Suppose that $Q_{2, 5}$ is a local jump. Since $Q_{2, 5}$ cannot be bad, we have that $N_{Q^*_{2, 5}}(c_4)\ne \emptyset$. Thus the shortest $c_2c_4$-path with interior in $Q^*_{2, 5}$ is a bad $v$-jump. If $Q_{2, 5}$ is not a local jump, then Lemma~\ref{lem-3-4}($Q_{2, 5}$) has a bad jump. Both contradict (\ref{eqa-3-3-1}). Thus we further have that ${\cal N}\subseteq (V(C)\setminus \{c_2, c_3, c_7\})\cup X_6\cup X_{1,6}\cup (X_{1,3}\setminus N(c_3))$.

Suppose that $[(X_{1,3}\setminus N(c_3))\cup (X_{1,6}\setminus N(c_6))\cup\{c_1\}]\cap {\cal N}\ne\emptyset$.  Let $Q_{1, 5}$ be a $c_1c_5$-path with shortest length and $Q^*_{1, 5}\subseteq {\cal D}\cup (X_{1,3}\setminus N(c_3))\cup (X_{1,6}\setminus N(c_6))$. Then $\{c_2, c_3, c_7\}$ is anticomplete to $Q^*_{1, 5}$. If $Q_{1, 5}$ is a local jump then Lemma~\ref{lem-3-2}($P_1, Q_{1, 5}$) has a bad jump. If $Q_{1, 5}$ is not local, then $N_{Q^*_{1, 5}}(c_4)\ne \emptyset$, and the shortest path from  ($c_4, \{c_1, c_6\}$)-path, with  $P^*\subseteq Q^*_{1, 5}$, is a bad jump. Both contradict (\ref{eqa-3-3-1}). Therefore, ${\cal N}\subseteq (V(C)\setminus \{c_1, c_2, c_3, c_7\})\cup X_6\cup X'_6$. This proves (\ref{eqa-3-9-0}).

Since $G$ has no $P_3$-cutsets, we have that ${\cal N}\cap (X_6\cup X'_6)\ne\emptyset$. If $c_4\in {\cal N}$, then there is a local $v$-jump $Q_{4, 6}$ across $c_5$ with interior in ${\cal D}$, and so Lemma~\ref{lem-3-1}($P_2, Q_{4, 6}$) has a bad jump, contradicting (\ref{eqa-3-3-1}).  Thus, ${\cal N}\subseteq (X_6'\cup X_6)\cup \{c_5, c_6\}$. Notice that every two vertices in $(X_6'\cup X_6)\cup\{c_5\}$ are joined by an induced path of length six or eight with interior in $((X_{1,6}\cup X_{2,6})\setminus N(c_6))\cup\{c_1,c_2,c_3,c_4\}$. We have that $(X_6'\cup X_6)\cup \{c_5,c_6\}$ is a parity star-cutset. This proves Claim~\ref{clm-3-5}.\qed

\medskip

To finish the proof of Lemma~\ref{theo-3-1}, we need to verify the case that $T$ is a short $e$-jump across $c_1c_2$, and $C$ has short $v$-jumps across $c_3$.

\begin{claim}\label{clm-3-6}
Suppose that  $T$ is a short $e$-jump across $c_1c_2$ and $C$ has a short $v$-jump across $c_3$. Then Lemma~\ref{theo-3-1} holds.
\end{claim}	
\pf Let $Q_{2, 4}$ be a short $v$-jump  across $c_3$. By Lemma~\ref{coro-2-2}, we have that $C$ has no short $v$-jumps across $c_1$, and no short $e$-jumps across $c_1c_7$.
Similar to the proofs of Claim~\ref{clm-3-4} and \ref{clm-3-5}, we can deduce that
${\cal X}=X_{2,4}\cup X_{3,7}\cup X_{1,3}$. Thus $c_5$ has no neighbor in ${\cal X}\cup \{c_1,c_2,c_3,c_7\}$, otherwise each vertex in $N_{{\cal X}}(c_5)$ provides us with a short jump starting from $c_5$

Since $d(c_5)\ge 3$, we choose ${\cal D}$ to be the vertex set of a maximal connected induced subgraph such that $N_{{\cal D}}(c_5)\ne \emptyset$ and ${\cal D}\cap (V(C)\cup {\cal X})=\emptyset$, and so (\ref{eqa-3-3-1}) holds. Let $X_4=N(c_4)\cap X_{2,4}$.
We claim that
\begin{equation}\label{eqa-3-14-0}
{\cal N}\subseteq X_4\cup \{c_4,c_5,c_6\}.
\end{equation}

Suppose that $[(X_{3, 7}\setminus N(c_7))\cup (X_{1, 3}\setminus N(c_1))\cup\{c_3\}]\cap {\cal N}\ne \emptyset$. Let $Q_{3, 5}$ be a $c_3c_5$-path with shortest length and  $Q^*_{3, 5}\subseteq {\cal D}\cup (X_{3, 7}\setminus N(c_7))\cup (X_{1, 3}\setminus N(c_1))$. If $Q_{3, 5}$ is a local jump, then Lemma~\ref{lem-3-2}($Q_{3, 5}, T$) has a bad jump. If $Q_{3, 5}$ is not a local jump, then Lemma~\ref{lem-3-4}($Q_{3, 5}$) has a bad jump. Both contradict  (\ref{eqa-3-3-1}). Thus ${\cal N}\subseteq (V(C)\setminus\{c_3\})\cup X_{2,4}\cup (X_{3,7}\cap N(c_7))\cup (X_{1, 3}\cap N(c_1))$.

Suppose that $[(X_{3,7}\cap N(c_7))\cup\{c_7\}]\cap {\cal N}\ne \emptyset$. Let $Q_{5, 7}$ be a $c_5c_7$-path with shortest length and $Q^*_{5, 7}\subseteq {\cal D}\cup (X_{3,7}\cap N(c_7))$. Then $N_{Q_{5, 7}^*}(c_3)=\emptyset$ as $c_3\not\in {\cal N}$. If $Q_{5, 7}$ is a local jump then Lemmas~\ref{lem-3-2}$(Q_{5, 7}, T)$ has a bad jump. If $Q_{5, 7}$ is not a local jump then Lemma~\ref{lem-3-4}($Q_{5, 7}$) has a bad jump. Both contradict (\ref{eqa-3-3-1}). Thus ${\cal N}\subseteq (V(C)\setminus\{c_3, c_7\})\cup X_{2,4}\cup (X_{1, 3}\cap N(c_1))$.

Suppose that $[(X_{2,4}\setminus N(c_4))\cup\{c_2\}]\cap {\cal N}\ne \emptyset$. Let $Q_{2, 5}$ be a $c_2c_5$-path with shortest length and  $Q^*_{2, 5}\subseteq {\cal D}\cup (X_{2,4}\setminus N(c_4))$. Then $N_{Q_{2, 5}^*}(c_3)=N_{Q_{2, 5}^*}(c_7)=\emptyset$ as $c_3, c_7\not\in {\cal N}$. If $Q_{2, 5}$ is a local jump, then $N_{Q_{2, 5}^*}(c_4)\ne\emptyset$ as $Q_{2, 5}$ cannot be short by (\ref{eqa-3-3-1}). This implies that the shortest $c_2c_4$-path, with interior in $Q_{2, 5}^*$, is a bad jump. If $Q_{2, 5}$ is not a local jump, then Lemma~\ref{lem-3-4}($Q_{2, 5}$) has a bad jump. Both contradict  (\ref{eqa-3-3-1}). Thus we have that ${\cal N}\subseteq (V(C)\setminus\{c_2, c_3, c_7\})\cup X_4\cup (X_{1, 3}\cap N(c_1))$.

Suppose that $[(X_{1, 3}\cap N(c_1))\cup\{c_1\}]\cap {\cal N}\ne \emptyset$. Let $Q_{1, 5}$ be a $c_1c_5$-path with shortest length and  $Q^*_{1, 5}\subseteq {\cal D}\cup (X_{1, 3}\cap N(c_1))$. Then $N_{Q_{1, 5}^*}(c_3)=N_{Q_{1, 5}^*}(c_7)=N_{Q_{1, 5}^*}(c_2)=\emptyset$.
If $Q_{1, 5}$ is a local jump, then Lemma~\ref{lem-3-2}($P_1, Q_{1, 5}$)has a bad jump. If $Q_{1, 5}$ is not local and $N_{Q_{1, 5}^*}(c_4)\ne\emptyset$, then Lemma~\ref{lem-3-4}($Q_{1, 5}$) has a bad jump. Both  contradict (\ref{eqa-3-3-1}). Therefore, ${\cal N}\subseteq (V(C)\setminus\{c_1, c_2, c_3, c_7\})\cup X_4$. This proves  (\ref{eqa-3-14-0}).

Since $G$ has no $P_3$-cutsets, we have that ${\cal N}\cap X_4\neq\emptyset$. If $c_6\in {\cal N}$, then there is a local $v$-jump $Q_{4, 6}$ across $c_5$ with interior in ${\cal D}$ such that Lemma~\ref{lem-3-1}($P_2, Q_{4, 6}$) has a bad jump, contradicting (\ref{eqa-3-3-1}).  Thus ${\cal N}\subseteq X_4\cup \{c_4,c_5\}$. Then every pair of distinct vertices in $X_4\cup\{c_5\}$ are joined by an induced path of length six or eight with interior in $\{c_1,c_2,c_6,c_7\}\cup X_{2,4}\setminus N(c_4)$, and so $X_4\cup\{c_4,c_5\}$ is a parity star-cutset. This proves Claim~\ref{clm-3-6}, and completes the proof of Lemma~\ref{theo-3-1}. \qed

\medskip

The proofs of the following lemmas take the same idea as that of above Lemma~\ref{theo-3-1}.

\begin{lemma}\label{theo-3-2}
Suppose that $C$ is of type $2$. Then $G$ admits a parity star-cutset.
\end{lemma}	
\pf Let $P_1$ be a local $v$ jump across $c_2$, and $P_2$ be a local $e$-jump across $c_6c_7$. Let $T$ be a short jump with $T^*\subseteq P_1^*\cup P_2^*$ in Lemma~\ref{lem-3-2}($P_1, P_2$), such that $T$ is a $v$-jump across $c_1$ or an $e$-jump across $c_1c_2$ or $c_1c_7$. By the definition of type 2, we have that
\begin{equation}\label{eqa-no-c4c7-v-jump}
\mbox{$C$ has no local $v$-jumps across $c_4$ or $c_7$,}
\end{equation}
and so
\begin{equation}\label{eqa-c6-empty}
\mbox{$N_{P_2^*}(c_6)=\emptyset$, and $T$ is not a jump across $c_1c_7$.}
\end{equation}

\begin{claim}\label{clm-3-7}
Suppose that $T$  is a short $v$-jump across $c_1$. Then Lemma~$\ref{theo-3-2}$ holds.
\end{claim}		
\pf It is certain that both $P_1$ and $P_2$ are not short. By the minimality of $|P_1^*\cup P_2^*|$, $C$ neither has short $v$-jumps across $c_2$ nor short $e$-jumps across $c_6c_7$. By Lemmas~\ref{coro-2-1} and \ref{theo-3-1}, we may assume that $C$ has no short $v$-jumps across any vertex in $\{c_3, c_4, c_5, c_6\}$. Thus except those across $c_1$, $C$ has no short $v$-jumps.

By Lemmas~\ref{coro-2-2}, we have that $C$ has no short $e$-jumps across $c_3c_4$ or $c_5c_6$. By Lemma~\ref{lem-3-3}, $C$ has no short jump across $c_2c_3$. Hence ${\cal X}=X_{2,7}\cup X_{2,6}\cup X_{3,7}$, and $c_5$ has no neighbor in ${\cal X}\cup \{c_1,c_2,c_3,c_7\}$ as $C$ has no short jumps starting from $c_5$.

Since $d(c_5)\ge 3$, we choose ${\cal D}$ to be the vertex set of a maximal connected induced subgraph such that $N_{{\cal D}}(c_5)\ne \emptyset$ and ${\cal D}\cap (V(C)\cup {\cal X})=\emptyset$. Let $X_6=X_{2,6}\cap N(c_6)$.
Now we prove that
\begin{equation}\label{eqa-3-19-0}
{\cal N}\subseteq X_6\cup \{c_4,c_5,c_6\}.
\end{equation}

Suppose that $[(X_{3,7}\cap N(c_3))\cup\{c_3\}]\cap {\cal N}\ne \emptyset$. Let $Q_{3, 5}$ be a $c_3c_5$-path with shortest length and $Q^*_{3, 5}\subseteq {\cal D}\cup (X_{3,7}\cap N(c_3))$. Since by (\ref{eqa-no-c4c7-v-jump}) $C$ has no local $v$-jump across $c_4$, we have that $N_{Q^*_{3, 5}}(\{c_1, c_2, c_6, c_7\})\ne\emptyset$, and so Lemma~\ref{lem-3-4}($Q_{3, 5}$) has a bad jump, contradicting  (\ref{eqa-3-3-1}). Thus ${\cal N}\subseteq (V(C)\setminus\{c_3\})\cup X_{2,7}\cup X_{2,6}\cup (X_{3,7}\setminus N(c_3))$.

Suppose that $[(X_{2,6}\setminus N(c_6))\cup (X_{2,7}\cap N(c_2))\cup\{c_2\}]\cap {\cal N}\ne \emptyset$. Let $Q_{2, 5}$ be a $c_2c_5$-path with shortest length and $Q^*_{2, 5}\subseteq {\cal D}\cup (X_{2,6}\setminus N(c_6))\cup (X_{2,7}\cap N(c_2))$. If $Q_{2, 5}$ is a local jump, then Lemma~\ref{lem-3-2}($Q_{2, 5}, T$) has a bad jump. If $Q_{2, 5}$ is not a local jump, then Lemma~\ref{lem-3-4}($Q_{2, 5}$) has a bad jump. Both contradict  (\ref{eqa-3-3-1}). Thus ${\cal N}\subseteq (V(C)\setminus\{c_2, c_3\})\cup X_6\cup (X_{2,7}\setminus N(c_2))\cup (X_{3,7}\setminus N(c_3))$.

Suppose that $c_1\in {\cal N}$. Let $Q_{1, 5}$ be a $c_1c_5$-path with shortest length and $Q^*_{1, 5}\subseteq {\cal D}\cup \{c_1\}$. Then $N_{Q^*_{1, 5}}(c_2)=N_{Q^*_{1, 5}}(c_3)=\emptyset$ as $c_2, c_3\not\in {\cal N}$. If $Q_{1, 5}$ is local, then Lemma~\ref{lem-3-2}($P_1, Q_{1, 5}$) has a bad jump. If $Q_{1, 5}$ is not local, then Lemma~\ref{lem-3-4}($Q_{1, 5}$) has a bad jump. Both contradict (\ref{eqa-3-3-1}). Thus ${\cal N}\subseteq (V(C)\setminus\{c_1, c_2, c_3\})\cup X_6\cup (X_{2,7}\setminus N(c_2))\cup (X_{3,7}\setminus N(c_3))$.

Suppose that $[(X_{2,7}\setminus N(c_2))\cup (X_{3,7}\setminus N(c_3))\cup\{c_7\}]\cap {\cal N}\ne \emptyset$. Let $Q_{5, 7}$ be a $c_5c_7$-path with shortest length and $Q^*_{5, 7}\subseteq {\cal D}\cup (X_{2,7}\setminus N(c_2))\cup (X_{3,7}\setminus N(c_3))$. Then $N_{Q^*_{5, 7}}(c_1)=N_{Q^*_{5, 7}}(c_2)=N_{Q^*_{5, 7}}(c_3)=\emptyset$. If $N_{Q^*_{5, 7}}(c_4)=\emptyset$, then let $Q'=Q_{5, 7}$. Otherwise, let $Q'$ be the shortest $c_4c_7$-path with interior in $Q^*_{5, 7}$. Then $Q'$ is a local jump, and so Lemma~\ref{lem-3-2}($Q', T$) has a bad jump, contradicting (\ref{eqa-3-3-1}). Therefore, ${\cal N}\subseteq (V(C)\setminus\{c_1, c_2, c_3, c_7\})\cup X_6$. This proves (\ref{eqa-3-19-0}).

Since $G$ has no $P_3$-cutsets, we have that ${\cal N}\cap X_6\ne\emptyset$. Then there is a short $e$-jump $Q_{2, 6}$ across $c_1c_7$. If $c_4\in {\cal N}$, then there is a local $v$-jump $Q_{4, 6}$ across $c_5$ with interior in ${\cal D}$, and Lemma~\ref{lem-3-2}($Q_{2, 6}, Q_{4, 6}$) has a bad jump. Hence ${\cal N}\subseteq X_6\cup \{c_5,c_6\}$. Since every two vertices in $X_6\cup\{c_5\}$ are joined by an induced path of length six or eight with interior in $X_{2,6}\setminus N(c_6)\cup\{c_2,c_3,c_4\}$, we have that $X_6\cup\{c_5,c_6\}$ is a parity star-cutset. This proves Claim~\ref{clm-3-7}.\qed

\medskip

By (\ref{eqa-c6-empty}), now suppose that $T$ is a short $e$-jump across $c_1c_2$. By Lemma~\ref{lem-3-2}, we have that
\begin{equation}\label{eqa-c4-local}
\mbox{$C$ has no local $v$-jumps across $c_4$ or $c_6$}.
\end{equation}

\begin{claim}\label{clm-3-8}
Suppose that $C$ has no short $v$-jump across $c_3$. Then Lemma~\ref{theo-3-2} holds.
\end{claim}	
\pf With the same arguments as that used in the proof of Claim~\ref{clm-3-7}, we have that no short jumps may contain $c_4$ or $c_5$. Thus
${\cal X}=X_{2,7}\cup X_{2,6}\cup X_{3,7}\cup X_{1,3}$, and $c_5$ is anticomplete to ${\cal X}\cup \{c_1,c_2,c_3,c_7\}$ as the vertices in $N_{{\cal X}}(c_5)$ may produce short jumps starting from $c_5$.

Since $d(c_5)\ge 3$, we choose ${\cal D}$ to be the set of a maximal connected induced subgraph such that $N_{{\cal D}}(c_5)\ne \emptyset$ and ${\cal D}\cap (V(C)\cup {\cal X})=\emptyset$. Then (\ref{eqa-3-3-1}) still holds. Let $X_6=X_{2,6}\cap N(c_6)$. We claim that
\begin{equation}\label{eqa-3-23-0}
{\cal N}\subseteq X_6\cup \{c_4,c_5,c_6\}.
\end{equation}

Suppose that $[(X_{3,7}\setminus N(c_7))\cup (X_{1,3}\cap N(c_3))\cup\{c_3\}]\cap {\cal N}\ne \emptyset$. Let $Q_{3, 5}$ be a $c_3c_5$-path with shortest length and $Q^*_{3, 5}\subseteq {\cal D}\cup (X_{3,7}\setminus N(c_7))\cup (X_{1,3}\cap N(c_3))$. Since by (\ref{eqa-c4-local}) $Q_{3, 5}$ is not a local jump, we have that Lemma~\ref{lem-3-4}($Q_{3, 5}$) has a bad jump, contradicting  (\ref{eqa-3-3-1}). Thus
${\cal N}\subseteq (V(C)\setminus\{c_3\})\cup X_{2,7}\cup X_{2,6}\cup (X_{3,7}\cap N(c_7))\cup (X_{1,3}\setminus N(c_3))$.

Suppose that $[(X_{2,7}\setminus N(c_2))\cup (X_{3,7}\cap N(c_7))\cup\{c_7\}]\cap {\cal N}\ne \emptyset$. Let $Q_{5, 7}$ be a $c_5c_7$-path with shortest length and $Q^*_{5, 7}\subseteq {\cal D}\cup (X_{2,7}\setminus N(c_2))\cup (X_{3,7}\cap N(c_7))$. Then $N_{Q_{5, 7}^*}(c_3)=\emptyset$. Since by (\ref{eqa-c4-local}) $Q_{5, 7}$ is not a local jump, thus Lemma~\ref{lem-3-4}($Q_{5, 7}$) has a bad jump, contradicting  (\ref{eqa-3-3-1}). So, ${\cal N}\subseteq (V(C)\setminus\{c_3, c_7\})\cup (X_{2,7}\cap N(c_2))\cup X_{2,6}\cup (X_{1,3}\setminus N(c_3))$.

Suppose that $[(X_{2,7}\cap N(c_2))\cup (X_{2,6}\setminus N(c_6))\cup\{c_2\}]\cap {\cal N}\ne \emptyset$. Let $Q_{2, 5}$ be a $c_2c_5$-path with shortest length and $Q^*_{2, 5}\subseteq {\cal D}\cup (X_{2,7}\cap N(c_2))\cup (X_{2,6}\setminus N(c_6))$. Then $N_{Q_{2, 5}^*}(c_3)=N_{Q_{2, 5}^*}(c_7)=\emptyset$. By Lemma~\ref{coro-2-2}, we see that $Q_{2, 5}$ is not a short jump.
If $Q_{2, 5}$ is a local jump, then $N_{Q_{2, 5}^*}(c_4)\ne\emptyset$, and the shortest $c_2c_4$-path with interior in $Q_{2, 5}^*$ is a bad jump, contradicting (\ref{eqa-3-3-1}). Thus $Q_{2, 5}$ is not a local jump, and so Lemma~\ref{lem-3-4}($Q_{2, 5}$) has a bad jump, contradicting (\ref{eqa-3-3-1}). Hence we have that  ${\cal N}\subseteq (V(C)\setminus\{c_2, c_3, c_7\})\cup X_{6}\cup (X_{1,3}\setminus N(c_3))$.

Suppose that $[(X_{1,3}\setminus N(c_3))\cup\{c_1\}]\cap {\cal N}\ne \emptyset$. Let $Q_{1, 5}$ be a $c_1c_5$-path with shortest length and $Q^*_{1, 5}\subseteq {\cal D}\cup (X_{1,3}\setminus N(c_3))$. Then $N_{Q_{1, 5}^*}(c_2)=N_{Q_{1, 5}^*}(c_3)=N_{Q_{1, 5}^*}(c_7)=\emptyset$. If $Q_{1, 5}$ is a local jump, then Lemma~\ref{lem-3-2}($P_1, Q_{1, 5}$) has a bad jump. If $Q_{1, 5}$ is not local, then Lemma~\ref{lem-3-4}($Q_{1, 5}$) has a bad jump. Both contradict (\ref{eqa-3-3-1}). Therefore, ${\cal N}\subseteq (V(C)\setminus\{c_1, c_2, c_3, c_7\})\cup X_{6}$. This proves (\ref{eqa-3-23-0}).

Since $G$ has no $P_3$-cutsets, we have that ${\cal N}\cap X_6\ne \emptyset$. Consequently, $X_{2, 6}\ne\emptyset$ and $C$ has a short jump, say $Q_{2,6}$ across $c_1c_7$. If $N_{\cal D}(c_4)\ne \emptyset$, then $C$ has a local jump $Q_{4, 6}$ across $c_5$ with interior in ${\cal D}$, and Lemma~\ref{lem-3-2}($Q_{2,6}, Q_{4, 6}$) has a bad jump, contradicting (\ref{eqa-3-3-1}).  Thus ${\cal N}\subseteq X_6\cup \{c_5,c_6\}$. Notice that each pair of vertices in $X_6\cup \{c_5\}$ are joined by an induced path of length six or eight with interior in $\{c_2,c_3,c_4\}\cup X_{2,6}\setminus N(c_6)$. Hence $X_6\cup \{c_5,c_6\}$ is a parity star-cutset. This proves Claim~\ref{clm-3-8}.\qed

\begin{claim}\label{clm-3-9}
Lemma~$\ref{theo-3-2}$ holds if  $C$ has a short $v$-jump across $c_3$.
\end{claim}	
\pf Suppose that $C$ has a short $v$-jump $Q_{2, 4}$ across $c_3$. By Lemmas~\ref{coro-2-2} and \ref{theo-3-1}, we have that
\begin{equation}\label{eqa-no-local-c5}
\mbox{$C$ has neither local $v$-jumps across $c_1$ or $c_5$, nor short $e$-jumps across $c_1c_7$ or $c_5c_6$.}
\end{equation}

With the similar arguments as that used in the proofs of Claims~\ref{clm-3-4} and \ref{clm-3-5}, we conclude that ${\cal X}=X_{2,4}\cup X_{3,7}\cup X_{1,3}$, and $c_5$ has no neighbors in ${\cal X}\cup \{c_1,c_2,c_3,c_7\}$.

Since $d(c_5)\ge 3$, we choose ${\cal D}$ to be the vertex set of a maximal connected induced subgraph such that $N_{{\cal D}}(c_5)\ne \emptyset$ and ${\cal D}\cap (V(C)\cup {\cal X})=\emptyset$.  Then (\ref{eqa-3-3-1}) still holds. Let $X_4=X_{2,4}\cap N(c_4)$.
We will prove that
\begin{equation}\label{eqa-3-28-0}
{\cal N}\subseteq X_4\cup \{c_4,c_5,c_6\}.
\end{equation}

Suppose that $[(X_{1,3}\setminus N(c_1))\cup (X_{3,7}\setminus N(c_7))\cup\{c_3\}]\cap {\cal N}\ne \emptyset$.  Let $Q_{3, 5}$ be a $c_3c_5$-path with shortest length and $Q^*_{3, 5}\subseteq {\cal D}\cup (X_{1,3}\setminus N(c_1))\cup (X_{3,7}\setminus N(c_7))$. Since $Q_{3, 5}$ is not a local jump  by (\ref{eqa-c4-local}), we have that Lemma~\ref{lem-3-4}($Q_{3, 5}$) has a bad jump, contradicting (\ref{eqa-3-3-1}). This shows that
${\cal N}\subseteq (V(C)\setminus\{c_3\})\cup X_{2,4}\cup (X_{3,7}\cap N(c_7))\cup (X_{1,3}\cap N(c_1))$.

Suppose that $[(X_{3,7}\cap N(c_7))\cup\{c_7\}]\cap {\cal N}\ne\emptyset$.  Let $Q_{5, 7}$ be a $c_5c_7$-path with shortest length and $Q^*_{5, 7}\subseteq {\cal D}\cup (X_{3,7}\cap N(c_7))$. Then $N_{Q_{5, 7}^*}(c_3)=\emptyset$.
Since by (\ref{eqa-c4-local}) $Q_{5, 7}$ is not a local jump, we have that Lemma~\ref{lem-3-4}($Q_{5, 7}$) has a bad jump, contradicting (\ref{eqa-3-3-1}). Thus ${\cal N}\subseteq (V(C)\setminus\{c_3, c_7\})\cup X_{2,4}\cup (X_{1,3}\cap N(c_1))$.

Suppose that $[(X_{2,4}\setminus N(c_4))\cup\{c_2\}]\cap {\cal N}\ne \emptyset$. Let $Q_{2, 5}$ be a $c_2c_5$-path with shortest length and $Q^*_{2, 5}\subseteq {\cal D}\cup (X_{2,4}\setminus N(c_4))$. Then $N_{Q_{2, 5}^*}(c_3)=N_{Q_{2, 5}^*}(c_7)=\emptyset$, and by (\ref{eqa-3-3-1}) $Q_{2, 5}$ is not a short jump. If $Q_{2, 5}$ is a local jump, then $N_{Q_{2, 5}^*}(c_4)\ne\emptyset$, and the shortest $c_2c_4$-path with interior in $Q_{2, 5}^*$ is a bad jump, contradicting (\ref{eqa-3-3-1}). Therefore, $Q_{2, 5}$ is not a local jump, and so Lemma~\ref{lem-3-4}($Q_{2, 5}$) has a bad jump, contradicting (\ref{eqa-3-3-1}). This shows that ${\cal N}\subseteq (V(C)\setminus\{c_2, c_3, c_7\})\cup X_{4}\cup (X_{1,3}\cap N(c_1))$.

Suppose that $[(X_{1,3}\cap N(c_1))\cup\{c_1\}]\cap {\cal N}\ne \emptyset$. Let $Q_{1, 5}$ be a $c_1c_5$-path with shortest length and $Q^*_{1, 5}\subseteq {\cal D}\cup (X_{1,3}\cap N(c_1))$. Then $N_{Q_{1, 5}^*}(c_2)=N_{Q_{1, 5}^*}(c_3)=N_{Q_{1, 5}^*}(c_7)=\emptyset$, and by (\ref{eqa-3-3-1}) $Q_{1, 5}$ is not a short jump. Now the shortest ($c_1, \{c_4, c_6\}$)-path, with interior in $Q_{1, 5}^*$, is a bad jump, contradicting (\ref{eqa-3-3-1}). Therefore, ${\cal N}\subseteq (V(C)\setminus\{c_1, c_2, c_3, c_7\})\cup X_{4}$, and (\ref{eqa-3-28-0}) holds.

Since $G$ has no $P_3$-cutsets, we have that $|{\cal N}\cap X_4|\ge 1$. Since $C$ has no local $v$-jump across $c_5$ by (\ref{eqa-no-local-c5}), we have that $N_{\cal D}(c_6)=\emptyset$, and so ${\cal N}\subseteq X_4\cup \{c_4,c_5\}$. Since each pair of vertices in $X_4\cup\{c_5\}$ are joined by an induced path of length six or eight with interior in $\{c_1,c_2,c_6,c_7\}\cup X_{2,4}\setminus N(c_4)$, we have that $X_4\cup\{c_4,c_5\}$ is a parity star-cutset. This proves Claim ~\ref{clm-3-9}, and also completes the proof of  Lemma~\ref{theo-3-2}. \qed

\begin{lemma}\label{theo-3-3}
Suppose that $C$ is of type $3$. Then $G$ admits a parity star-cutset.
\end{lemma}	
\pf Let $P_1$ be a local $e$-jump across $c_2c_3$, and $P_2$ be a local $e$-jump across $c_6c_7$. It follows from the definition of type 3, $C$ has no local $v$-jumps across any vertex in $\{c_2, c_4, c_5, c_7\}$, and so $N_{P_1^*}(c_3)=N_{P_2^*}(c_6)=\emptyset$. Consequently, Lemma~\ref{lem-3-3}($P_1, P_2$) cannot have short jumps across $c_1c_2$ or $c_1c_7$.  Let $T$ be a short $v$-jump across $c_1$ in Lemma~\ref{lem-3-3}($P_1, P_2$). Then none of $P_1$ and $P_2$ can be short by Lemma~\ref{lem-3-3}.

By Lemma~\ref{coro-2-1}, we have that $C$ has no short $v$-jumps across $c_3$ or $c_6$. By Lemma~\ref{theo-3-2}, we may suppose that $C$ has no  local $e$-jumps across $c_3c_4$ or $c_5c_6$. Therefore,
\begin{equation}\label{eqa-v-c1-jump}
\mbox{the only possible local $v$-jumps are those across $c_1$,}
\end{equation}
and
\begin{equation}\label{eqa-e-c1c2-jump}
\mbox{the only possible local $e$-jumps are those across some edge in $\{c_1c_2, c_2c_3, c_6c_7\}$}.
\end{equation}
Consequently, ${\cal X}=X_{2,7}\cup X_{3,7}\cup X_{2,6}$, and $c_5$ has no neighbor in ${\cal X}\cup \{c_1,c_2,c_3,c_7\}$ as each vertex in $N_{{\cal X}}(c_5)$ provides us with a short jump starting from $c_5$.

Since $d(c_5)\ge 3$, we  choose ${\cal D}$ to be the vertex set of a maximal connected induced subgraph such that $N_{{\cal D}}(c_5)\ne \emptyset$ and ${\cal D}\cap (V(C)\cup {\cal X})=\emptyset$. Then, (\ref{eqa-3-3-1}) still holds. Let $X_6=X_{2, 6}\cap N(c_6)$. We will prove that
\begin{equation}\label{eqa-3-33-0}
{\cal N}\subseteq X_6\cup \{c_4,c_5,c_6\}.
\end{equation}

Suppose that $[(X_{3,7}\cap N(c_3))\cup\{c_3\}]\cap {\cal N}\ne \emptyset$, and let $Q_{3, 5}$ be a $c_3c_5$-path with shortest length and $Q^*_{3, 5}\subseteq {\cal D}\cup (X_{3,7}\cap N(c_3))$. Since $C$ has no local $v$-jump across $c_4$ by (\ref{eqa-v-c1-jump}), we have that Lemma~\ref{lem-3-4}($Q_{3, 5}$) has a bad jump, contradicting (\ref{eqa-3-3-1}). Thus ${\cal N}\subseteq (V(C)\setminus\{c_3\})\cup X_{2,7}\cup X_{2,6}\cup (X_{3,7}\setminus N(c_3))$.

Suppose that $[(X_{2,6}\setminus N(c_6))\cup (X_{2,7}\cap N(c_2))\cup\{c_2\}]\cap {\cal N}\ne \emptyset$,  and let $Q_{2, 5}$ be a $c_2c_5$-path with shortest length and $Q^*_{2, 5}\subseteq {\cal D}\cup (X_{2,6}\setminus N(c_6))\cup (X_{2,7}\cap N(c_2))$. Since $Q_{2, 5}$ is not a local jump by (\ref{eqa-e-c1c2-jump}), Lemma~\ref{lem-3-4}($Q_{2, 5}$) has a bad jump, contradicting (\ref{eqa-3-3-1}). Consequently, ${\cal N}\subseteq (V(C)\setminus\{c_2, c_3\}) \cup X_{6}\cup (X_{2,7}\setminus N(c_2))\cup (X_{3,7}\setminus N(c_3))$.

Suppose $c_1\in {\cal N}$. Let $Q_{1, 5}$ be a $c_1c_5$-path with shortest length and $Q^*_{1, 5}\subseteq {\cal D}$. If $Q_{1, 5}$ is a local jump then Lemma~\ref{lem-3-3}($P_1, Q_{1, 5}$) has a bad jump. If $Q_{1, 5}$ is not a local jump then Lemma~\ref{lem-3-4}($Q_{1, 5}$) has a bad jump. Both contradict (\ref{eqa-3-3-1}). Thus ${\cal N}\subseteq (V(C)\setminus\{c_1, c_2, c_3\}) \cup X_{6}\cup (X_{2,7}\setminus N(c_2))\cup (X_{3,7}\setminus N(c_3))$.

Suppose that $[(X_{2,7}\setminus N(c_2))\cup (X_{3,7}\setminus N(c_3))\cup\{c_7\}]\cap {\cal N}\ne \emptyset$, and let $Q_{5, 7}$ be a $c_5c_7$-path with shortest length and $Q^*_{5, 7}\subseteq {\cal D}\cup (X_{2,7}\setminus N(c_2))\cup (X_{3,7}\setminus N(c_3))$. Then $N_{Q^*_{5, 7}}(c_1)=N_{Q^*_{5, 7}}(c_2)=N_{Q^*_{5, 7}}(c_3)=\emptyset$. If $N_{Q^*_{5, 7}}(c_4)=\emptyset$, then $Q_{5, 7}$ is a local $v$-jump, contradicting (\ref{eqa-v-c1-jump}). Otherwise,  $N_{Q^*_{5, 7}}(c_4)\ne \emptyset$ and $C$ has a local $e$-jump across $c_5c_6$, contradicting (\ref{eqa-e-c1c2-jump}). Therefore, ${\cal N}\subseteq (V(C)\setminus\{c_1, c_2, c_3, c_7\}) \cup X_{6}$. This proves (\ref{eqa-3-33-0}).

Since $G$ has no $P_3$-cutsets, we have that ${\cal N}\cap X_6\ne\emptyset$. If $c_4\in {\cal N}$, then $C$ has a local $v$-jump $Q_{4, 6}$ across $c_5$, contradicting (\ref{eqa-v-c1-jump}). Thus ${\cal N}\subseteq X_6\cup \{c_5,c_6\}$. Since each pair of vertices in $X_6\cup\{c_5\}$ are joined by an induced path of length six or eight with interior in $\{c_2,c_3,c_4\}\cup X_{2,6}\setminus N(c_6)$, we have that $X_6\cup\{c_5,c_6\}$ is a parity star-cutset. This proves Lemma~\ref{theo-3-3}. \qed

\medskip

Finally, we consider the case where $C$ has at least two kinds of equivalent local jumps, and is not of  type $i$ for any $i\in \{1, 2, 3\}$.

\begin{lemma}\label{theo-3-4}
Suppose that $C$ is not of type $i$ for any $i\in \{1, 2, 3\}$, and $C$ has at least two kinds of equivalent local jumps. If $C$ has a local $v$-jump, then $G$ admits a parity star-cutset.
\end{lemma}	
\pf Without loss of generality, suppose that $C$ has a local $v$-jump across $c_3$. Let $P_1$ be a shortest local jump across $c_3$. By Lemmas~\ref{theo-3-1} and \ref{theo-3-2}, we may assume that
\begin{equation}\label{eqa-3-4-local-jumps}
\mbox{$C$ has no local $v$-jumps across $c_1$ or $c_5$, and no local $e$-jumps across $c_1c_7$ or $c_5c_6$.}
\end{equation}

By symmetry, we need to consider the situations that $C$ has a local $v$-jump across one vertex in $\{c_2, c_7\}$, or a local $e$-jump across one edge in $\{c_1c_2, c_3c_4, c_6c_7\}$.

\begin{claim}\label{clm-4-1}
$C$ has no local $v$-jumps across $c_7$.
\end{claim}	
\pf Suppose to its contrary, let $P_2$ be a local $v$-jump across $c_7$ shortest length. By Lemmas~\ref{theo-3-1} and \ref{theo-3-2}, suppose that
\begin{equation}\label{eqa-2-3-local-jumps}
\mbox{$C$ has no local $v$-jumps across $c_2$, and no local $e$-jumps across $c_2c_3$ or $c_4c_5$.}
\end{equation}

Firstly we prove that 	
\begin{equation}\label{eqa-4-1-1}
\mbox{$C$ has a short $e$-jump $Q_{3, 7}$ across $c_1c_2$, with interior in $P_1^*\cup P_2^*$.}
\end{equation}		

Since $G$ induces no big odd holes, we have that $P_1^*$ cannot be  disjoint from and anticomplete to $P_2^*$.  Let $Q_{2, 6}$ be the shortest $c_2c_6$-path with interior in $P_1^*\cup P_2^*$. Then $Q_{2, 6}$ is not a local jump by (\ref{eqa-3-4-local-jumps}). If $N_{Q^*_{2, 6}}(c_4)\ne\emptyset$, then $C$ has a local $v$-jump across $c_5$ or a local $e$-jump across $c_2c_3$ or $c_5c_6$, contradicting (\ref{eqa-3-4-local-jumps}) or (\ref{eqa-2-3-local-jumps}). Thus $N_{Q^*_{2, 6}}(c_3)\ne\emptyset$. Let $Q_{3,6}$ be the shortest $c_3c_6$-path with interior in $Q^*_{2, 6}$. By (\ref{eqa-2-3-local-jumps}), we have that $Q_{3,6}$ cannot be a local jump, which implies that  $N_{Q_{3,6}}(c_7)\ne\emptyset$. This implies that $C$ has a short $e$-jump across $c_1c_2$, with interior in $P_1^*\cup P_2^*$. Therefore, (\ref{eqa-4-1-1}) holds.

Consequently, we have that neither $P_1$ nor $P_2$ can be short, that is
\begin{equation}\label{eqa-3-7-local-jumps}
\mbox{$C$ has no local $v$-jumps across $c_3$ or $c_7$.}
\end{equation}
By Lemmas~\ref{coro-2-2}, \ref{coro-2-3}, \ref{theo-3-2}, and \ref{theo-3-3}, we may assume that
\begin{equation}\label{eqa-4-5-6-local-jumps}
\mbox{$C$ has no local $v$-jumps across $c_4$ or $c_6$.}
\end{equation}

If $C$ has a short $e$-jump $Q_{1, 5}$ across $c_6c_7$, then $Q^*_{1, 5}$ is disjoint and anticomplete to $Q^*_{3, 7}$, and so $c_1c_7Q_{3, 7}c_3c_4c_5Q_{1, 5}$ is a big odd hole, a contradiction. By symmetry, we may assume that
\begin{equation}\label{eqa-34-67-short}
\mbox{$C$ has no short $e$-jumps across $c_3c_4$ or $c_6c_7$.}
\end{equation}
Combining this with (\ref{eqa-3-4-local-jumps}), (\ref{eqa-2-3-local-jumps}), (\ref{eqa-3-7-local-jumps}) and (\ref{eqa-4-5-6-local-jumps}), we have that $C$ has no short $v$-jumps, and the only possible short $e$-jumps are those across $c_1c_2$.  Thus, ${\cal X}=X_{3,7}$, and $c_5$ has no neighbor in ${\cal X}\cup \{c_1,c_2,c_3,c_7\}$ as each vertex in $N_{{\cal X}}(c_5)$ provides us with a short jump starting from $c_5$.

Since $d(c_5)\ge 3$, we choose ${\cal D}$ to be the vertex set of a maximal connected induced subgraph such that $N_{{\cal D}}(c_5)\ne \emptyset$ and ${\cal D}\cap (V(C)\cup {\cal X})=\emptyset$.  Then (\ref{eqa-3-3-1}) holds.

suppose that $[(X_{3,7}\cap N(c_3))\cup\{c_3\}]\cap {\cal N}\ne \emptyset$, and let $Q_{3, 5}$ be a $c_3c_5$-path with shortest length and $Q^*_{3, 5}\subseteq {\cal D}\cup (X_{3,7}\cap N(c_3))$. Since $Q_{3, 5}$ is not a local jump by (\ref{eqa-4-5-6-local-jumps}), Lemma~\ref{lem-3-4}($Q_{3, 5}$) has a bad jump, contradicting (\ref{eqa-3-3-1}). Therefore, ${\cal N}\subseteq (V(C)\setminus\{c_3\})\cup (X_{3,7}\setminus N(c_3))$.

If $c_2\in {\cal N}$, we can get a contradiction with the same argument as above. Thus $c_2\notin {\cal N}$. By symmetry, $[(X_{3,7}\setminus N(c_3))\cup\{c_7\}]\cap {\cal N}=\emptyset$ and $c_1\notin {\cal N}$. Therefore, ${\cal N}\subseteq \{c_4,c_5,c_6\}$, which leads to a contradiction to the choice of $G$. This proves Claim~\ref{clm-4-1}. \qed

\begin{claim}\label{clm-4-2}
Suppose that $C$ has a local $v$-jump across $c_2$. Then Lemma~$\ref{theo-3-4}$ holds.
\end{claim}	
\pf  We choose $P_2$ to be a local jump across $c_2$ with shortest length.  By Lemmas~\ref{theo-3-1} and  \ref{theo-3-2},  and by Claim~\ref{clm-4-1}, we may assume that
\begin{equation}\label{eqa-4-2-1}
\mbox{$C$ has no local $v$-jumps across $c_4$ or $c_6$ or $c_7$, and no local $e$-jumps across $c_4c_5$ or $c_6c_7$.}
\end{equation}

Since $G$ is big odd hole free, we have that $C$ cannot have both short $e$-jumps across $c_1c_2$ and short $e$-jumps across $c_3c_4$. Thus we may suppose, by symmetry, that
\begin{equation}\label{eqa-4-2-2}
\mbox{$C$ has no short $e$-jumps across $c_3c_4$.}
\end{equation}

By (\ref{eqa-3-4-local-jumps}) and (\ref{eqa-4-2-1}),  we have that
\begin{equation}\label{eqa-5-6-localjumps}
\mbox{$C$ has no local jumps with end either $c_5$ or $c_6$.}
\end{equation}
Thus the only possible short $v$-jumps of $C$ are those across $c_2$ or $c_3$, and only possible short $e$-jumps of $C$ are those across $c_2c_3$ or $c_1c_2$.  Hence ${\cal X}=X_{1,3}\cup X_{2,4}\cup X_{1,4}\cup X_{3,7}$, and both $c_5$ and $c_6$ have no neighbors in ${\cal X}\cup \{c_1,c_2,c_3,c_4\}$.

Since $d(c_6)\ge 3$, we choose ${\cal D}$ to be the vertex set of a maximal connected induced subgraph such that $N_{{\cal D}}(c_6)\ne \emptyset$ and ${\cal D}\cap (V(C)\cup {\cal X})=\emptyset$.   Then (\ref{eqa-3-3-1}) still holds. Let $X_7=X_{3,7}\cap N(c_7)$. We will prove that
\begin{equation}\label{eqa-4-2-3-0}
{\cal N}\subseteq X_7\cup \{c_5,c_6,c_7\}.
\end{equation}

Suppose that $[((X_{1,4}\cup X_{2,4})\cap N(c_4))\cup\{c_4\}]\cap {\cal N}\ne \emptyset$, and let $Q_{4, 6}$ be a  $c_4c_6$-path with shortest length and $Q^*_{4, 6}\subseteq {\cal D}\cup ((X_{1,4}\cup X_{2,4})\cap N(c_4))$. Since $Q_{4, 6}$ is not a local jump by (\ref{eqa-3-4-local-jumps}), Lemma~\ref{lem-3-4}($Q_{4, 6}$) has a bad jump, contradicting (\ref{eqa-3-3-1}). Hence
${\cal N}\subseteq (V(C)\setminus\{c_4\})\cup X_{1,3}\cup X_{3,7}\cup (X_{2,4}\cup X_{1,4})\setminus N(c_4)$.

Suppose that $[(X_{1,3}\cap N(c_3))\cup (X_{3,7}\setminus N(c_7))\cup\{c_3\}]\cap {\cal N}\ne\emptyset$, and let $P_{3, 6}$ be a $c_3c_6$-path with shortest length and $P^*_{3, 6}\subseteq {\cal D}\cup (X_{1,3}\cap N(c_3))\cup (X_{3,7}\setminus N(c_7))$. Since $P_{3, 6}$ is not a local jump by (\ref{eqa-4-2-1}), Lemma~\ref{lem-3-4}($Q_{3, 6}$) has a bad jump, contradicting (\ref{eqa-3-3-1}). Consequently,
${\cal N}\subseteq (V(C)\setminus\{c_3, c_4\})\cup X_{7}\cup (X_{1,3}\setminus N(c_3))\cup (X_{2,4}\cup X_{1,4})\setminus N(c_4)$.

suppose that $[(X_{2,4}\setminus N(c_4))\cup\{c_2\}]\cap {\cal N}\ne \emptyset$, and let $Q_{2, 6}$ be a $c_2c_6$-path with shortest length and $Q^*_{2, 6}\subseteq {\cal D}\cup (X_{2,4}\setminus N(c_4))$. Since $Q_{2, 6}$ is not a local jump by (\ref{eqa-5-6-localjumps}), Lemma~\ref{lem-3-4}($Q_{2, 6}$) has a bad jump, a contradiction to (\ref{eqa-3-3-1}).  Thus ${\cal N}\subseteq (V(C)\setminus\{c_2, c_3, c_4\})\cup X_{7}\cup (X_{1,3}\setminus N(c_3))\cup (X_{1,4}\setminus N(c_4))$.

Suppose that $[(X_{1,3}\setminus N(c_3))\cup (X_{1,4}\setminus N(c_4))\cup\{c_1\}]\cap {\cal N}\ne \emptyset$, and let $Q_{1, 6}$ be a $c_1c_6$-path with shortest length and $Q^*_{1, 6}\subseteq {\cal D}\cup (X_{1,3}\setminus N(c_3))\cup (X_{1,4}\setminus N(c_4))$. Since $Q_{1, 5}$ is not a local jump by (\ref{eqa-5-6-localjumps}), Lemma~\ref{lem-3-4}($Q_{1, 6}$) has a bad jump, which contradicts (\ref{eqa-3-3-1}). Therefore, ${\cal N}\subseteq (V(C)\setminus\{c_1, c_2, c_3, c_4\})\cup X_{7}$. This proves (\ref{eqa-4-2-3-0}).

Since $G$ has no $P_3$-cutsets, we have that $X_7\ne\emptyset$. Let $Q_{3, 7}$ be a short jump across $c_1c_2$. If $c_5\in {\cal N}$, then $C$ has a local jump $Q_{5, 7}$ across $c_6$ with interior in ${\cal D}$, and Lemma~\ref{lem-3-2}($Q_{3,7}, Q_{5, 7}$) has a bad jump, contradicting (\ref{eqa-3-3-1}).
Therefore, ${\cal N}\subseteq X_7\cup \{c_7, c_6\}$. Notice that each pair of vertices in $X_7\cup \{c_6\}$ are joined by an induced path of length six or eight with interior in $\{c_3,c_4,c_5\}\cup X_{3,7}\setminus N(c_7)$. We have that $X_7\cup \{c_6,c_7\}$ is a parity star-cutset. This proves Claim~\ref{clm-4-2}.\qed

\begin{claim}\label{clm-4-3}
Suppose that $C$ has a local $e$-jump across $c_6c_7$. Then Lemma~$\ref{theo-3-4}$ holds.
\end{claim}	
\pf We choose $P_2$ to be a local jump across $c_6c_7$ with shortest length.  By Claim~\ref{clm-4-1}, we may assume that $N_{P_2^*}(c_6)=N_{P_2^*}(c_7)=\emptyset$. Thus $P_2$ is short. Since $G$ is big odd hole free, we have that $P_1^*$ is not anticomplete to $P_2^*$. Then the shortest $c_2c_5$-path $Q_{2, 5}$, with interior in $P_1^*\cup P_2^*$, is a local jump, which together with $P_2$ implies that $C$ is of type 3, a contradiction. This proves Claim~\ref{clm-4-3}.\qed

\begin{claim}\label{clm-4-4}
Suppose that $C$ has a local $e$-jump across $c_1c_2$. Then Lemma~$\ref{theo-3-4}$ holds.
\end{claim}	
\pf We choose $P_2$ to be a local jump across $c_1c_2$ with shortest length.  By Lemmas~\ref{theo-3-1}, \ref{theo-3-2} and \ref{theo-3-3},  and by Claims~\ref{clm-4-1}, \ref{clm-4-2} and \ref{clm-4-3}, we may assume that $C$ has no local $v$-jumps across any vertex in $V(C)\setminus\{c_3\}$, and
no local $e$-jumps across any edge in $\{c_1c_7, c_4c_5, c_5c_6, c_6c_7\}$.

Since $G$ is big odd hole free, we have that $C$ cannot have both a short $e$-jump across $c_1c_2$ and a short $e$-jump across $c_3c_4$.

If $C$ has no short $e$-jumps across $c_3c_4$, then the only possible short $e$-jumps of $C$ are those across $c_1c_2$ or $c_2c_3$. Thus  ${\cal X}=X_{1,4}\cup X_{2,4}\cup X_{3,7}$, and $c_6$ has no neighbor in ${\cal X}\cup \{c_1,c_2,c_3,c_4\}$ as otherwise each vertex in $N_{{\cal X}}(c_6)$ provides us with a short jump starting from $c_6$.

If $C$ has no short $e$-jumps across $c_1c_2$, then the only possible short $e$-jumps are those across $c_3c_4$ or $c_2c_3$.  Thus ${\cal X}=X_{1,4}\cup X_{2,4}\cup X_{2,5}$, and $c_6$ has no neighbor in ${\cal X}\cup \{c_1,c_2,c_3,c_4\}$ as otherwise each vertex in $N_{{\cal X}}(c_6)$ provides us with a short jump starting from $c_6$.

In both cases above, we choose ${\cal D}$ to be the vertex set of a maximal connected induced subgraph with $N_{{\cal D}}(c_6)\ne \emptyset$ and ${\cal D}\cap (V(C)\cup {\cal X})=\emptyset$. Then (\ref{eqa-3-3-1}) still holds.

With a similar  argument as that used in the proof of Claim \ref{clm-4-2}, we conclude that
\begin{itemize}
\item  $(X_{3,7}\cap N(c_7))\cup \{c_6,c_7\}$ is a parity star-cutset if $C$ has no short $e$-jumps across $c_3c_4$, and

\item $(X_{2,5}\cap N(c_5))\cup \{c_5,c_6\}$ is a parity star-cutset if $C$ has no short $e$-jumps across $c_1c_2$.
\end{itemize}
This completes the proof of Claim~\ref{clm-4-4}.\qed

\begin{claim}\label{clm-4-5}
Suppose that $C$ has a local $e$-jump across $c_3c_4$. Then Lemma~$\ref{theo-3-4}$ holds.
\end{claim}	
\pf We choose $P_2$ to be a local $e$-jump across $c_3c_4$ with shortest length.  By Lemmas~\ref{theo-3-1}, \ref{theo-3-2} and \ref{theo-3-3},  and by Claim~\ref{clm-4-1}, \ref{clm-4-2}, \ref{clm-4-3} and  \ref{clm-4-4}, we may assume that the only possible local jumps of $C$ are those across $c_3$ or $c_2c_3$ or $c_3c_4$. Thus ${\cal X}=X_{1,4}\cup X_{2,4}\cup X_{2,5}$, and $c_6$ is anticomplete to ${\cal X}\cup \{c_1,c_2,c_3,c_4\}$ as otherwise any vertex in $N_{{\cal X}}(c_6)$ provides us with a local jump starting from $c_6$.

We choose ${\cal D}$ to be the vertex set of a maximal connected induced subgraph such that $N_{{\cal D}}(c_6)\ne \emptyset$ and ${\cal D}\cap (V(C)\cup {\cal X})=\emptyset$. Then (\ref{eqa-3-3-1}) still holds.

With a similar  argument as that used in the proof of Claim \ref{clm-4-2}, we conclude that
$(X_{2,5}\cap N(c_5))\cup \{c_5,c_6\}$ is a parity star-cutset. This proves Claim~\ref{clm-4-5}, and completes the proof of Lemma~\ref{theo-3-4}.\qed

\begin{lemma}\label{theo-3-5}
Suppose that $C$ is not of  type $i$ for any $i\in\{1, 2, 3\}$, and $C$ has at least two kinds of equivalent local jumps. If $C$ has no local $v$-jumps, then $G$ admits a parity star-cutset.
\end{lemma}	
\pf Suppose that $C$ has no local $v$-jumps. Without loss of generality, we suppose by Lemma~\ref{lem-2-1} that
\begin{equation}\label{eqa-no-local-v-local-34}
\mbox{$C$ has a local $e$-jump across $c_3c_4$,}
\end{equation}
and let $P_1$ be a local $e$-jump across $c_3c_4$ with shortest length.

Notice that $P_1$ is a local $e$-jump across $c_3c_4$ with shortest length, by Lemma~\ref{theo-3-3}, we may assume that
\begin{equation}\label{eqa-12-23-local-e-jumps}
\mbox{$C$ has no local $e$-jumps across $c_1c_7$ or $c_6c_7$.}
\end{equation}

By (\ref{eqa-no-local-v-local-34}) and (\ref{eqa-12-23-local-e-jumps}), and by symmetry, we need to consider the cases where $C$ has a local $e$-jump across $c_1c_2$ or $c_2c_3$.

\begin{claim}\label{clm-5-1}
Suppose that $C$ has a local $e$-jump across $c_1c_2$. Then Lemma~$\ref{theo-3-5}$ holds.
\end{claim}	
\pf We choose $P_2$ to be a local jump across $c_1c_2$ with shortest length.  By Lemmas~\ref{theo-3-1}, \ref{theo-3-2}, \ref{theo-3-3} and \ref{theo-3-4}, we may assume that
the only possible local jumps of $C$ are those across $c_1c_2$ or $c_2c_3$ or $c_3c_4$.

Since $G$ is big odd hole free, at most one of $P_1$ and $P_2$ is short. By symmetry, we suppose that $P_2$ is not a short jump. Then ${\cal X}=X_{1, 4}\cup X_{2, 5}$, and $c_6$ is anticomplete to ${\cal X}\cup \{c_1,c_2,c_3,c_4\}$ as otherwise any vertex in $N_{{\cal X}}(c_6)$ provides us with a local jump starting from $c_6$.

We choose ${\cal D}$ to be the vertex set of a maximal connected induced subgraph such that $N_{{\cal D}}(c_6)\ne \emptyset$ and ${\cal D}\cap (V(C)\cup {\cal X})=\emptyset$. Then (\ref{eqa-3-3-1}) still holds.

With a similar  argument as that used in the proof of Claim~\ref{clm-4-2}, we conclude that
$(X_{2,5}\cap N(c_5))\cup \{c_5,c_6\}$ is a parity star-cutset. This proves Claim~\ref{clm-5-1}.\qed

\begin{claim}\label{clm-5-2}
Suppose that $C$ has a local $e$-jump across $c_2c_3$. Then Lemma~$\ref{theo-3-5}$ holds.
\end{claim}	
\pf We choose $P_2$ to be a local jump across $c_2c_3$ with shortest length.  By Lemmas~\ref{theo-3-1}, \ref{theo-3-2}, \ref{theo-3-3} and \ref{theo-3-4}, and by Claim~\ref{clm-5-1}, we may assume that the only possible local jumps are those across $c_2c_3$ or $c_3c_4$.
Again we have that ${\cal X}=X_{1, 4}\cup X_{2, 5}$, and $c_6$ is anticomplete to ${\cal X}\cup \{c_1,c_2,c_3,c_4\}$.

We choose ${\cal D}$ to be the vertex set of a maximal connected induced subgraph such that $N_{{\cal D}}(c_6)\ne \emptyset$ and ${\cal D}\cap (V(C)\cup {\cal X})=\emptyset$. Then (\ref{eqa-3-3-1}) still holds.

With a similar  argument as that used in the proof of Claim~\ref{clm-4-2}, we conclude that
$(X_{2,5}\cap N(c_5))\cup \{c_5,c_6\}$ is a parity star-cutset. This proves Claim~\ref{clm-5-2}. \qed

\bigskip

\noindent{\bf Proof of Theorem~\ref{theo-1-3}}. Suppose to its contrary. By the conclusion of \cite{WXX2022}, we choose $G$ to be a minimal heptagraph with $\chi(G)=4$. Then, $\delta(G)\ge 3$, and $G$ is not bipartite. Let $C$ be a 7-hole of $G$. By Lemma~\ref{lem-critical-H}, $G$ has no $P_3$-cutsets or parity-star cutsets. By Lemmas~\ref{theo-1-4} and \ref{theo-1-5}, $G$ induces no ${\cal P}$ or ${\cal P}'$.
By Lemma~\ref{lem-unique-local}, $C$ has at least two kinds of equivalent local jumps.  By Lemmas~\ref{theo-3-1}, \ref{theo-3-2} and \ref{theo-3-3}, $C$ cannot be type $i$ for any $i\in\{1, 2, 3\}$. By Lemma~\ref{theo-3-4}, $C$ has no local $v$-jumps. This indicates that the only possible local jumps of $C$ must be $e$-jumps, which leads to a contradiction to Lemma~\ref{theo-3-5}. \qed

\bigskip

\noindent{\bf Remark.} Recall that ${\cal G}_{\l}$ is the family of graphs without cycles of length at most $2\l$ and without odd holes of length at least $2\l+3$. The current authors proposed a conjecture claiming that all graphs in $\cup_{\l\ge 2} {\cal G}_{\l}$ are 3-colorable.

It seems that the structures of graphs in ${\cal G}_{\l}$ have some connection with cages. For given integers $k$ and $g$, a $(k, g)$-{\em cage} is a $k$-regular graph which has girth $g$ and the smallest number of vertices. The unique $(3, 5)$-cage is the Petersen graph, and the unique $(3, 7)$-cage is the McGee graph \cite{WM1960, WT1966}. Notice that the graph ${\cal P}'$ can be obtained from the McGee graph by deleting four disjoint groups of vertices such that each group induces a path of length 2. The graph, obtained from the Petersen graph by deleting two adjacent vertices, plays an important role in \cite{MCPS2022}, and the graph ${\cal P}'$ is also crucial in the proof of Theorem~\ref{theo-1-3}.

Since the Balaban graph \cite{ATB1973, MMN1998} is the unique $(3, 11)$-cage with 112 vertices, perhaps one can prove that all graphs in ${\cal G}_{5}$ are 3-colorable following the idea of  Chudnovsky and Seymour \cite{MCPS2022},  with much more detailed analysis. But it seems very hard to solve the 3-colorability of graphs in ${\cal G}_{4}$ along this approach, as there are eighteen $(3, 9)$-cages each of which  has 58 vertices (see \cite{GERJ2013, PW1982}).

Nelson, Plummer, Robertson, and Zha, \cite{NPRZ2011} proved that the Petersen graph is the only non-bipartite cubic pentagraph which is 3-connected and internally 4-connected, and Plummer and Zha \cite{MPXZ} presented some  3-connected and internally 4-connected non-bipartite non-cubic pentagraphs. It is known that for $g\ge 5$, all $(3, g)$-cages are 3-connected and internally 4-connected (see \cite{MPB2003}). Notice that the McGee graph has 9-holes, and so is not in ${\cal G}_{3}$. It seems interesting to consider the existences of non-bipartite 3-connected, internally 4-connected graphs for ${\cal G}_{\l}$ ($\l\ge 3$).

\end{document}